\documentclass[notitlepage, 11pt]{article}
\usepackage{graphicx}

\usepackage{epstopdf}
\usepackage{color}
\usepackage{latexsym}
\usepackage{amsmath, amssymb, amsfonts}
\usepackage[toc]{appendix}
\numberwithin{equation}{section}
\DeclareMathOperator*{\argmin}{arg\,min}
\DeclareMathOperator*{\argmax}{arg\,max}

\newtheorem{lemma}{Lemma}
\newtheorem{theorem}{Theorem}
\newtheorem{corollary}{Corollary}
\newtheorem{definition}{Definition}

\newtheorem{proposition}{Proposition}

\newtheorem{remark}{Remark}

\newcommand{\beginsec}{
\setcounter{lemma}{0}
\setcounter{theorem}{0}
\setcounter{corollary}{0}
\setcounter{definition}{0}
\setcounter{example}{0}
\setcounter{proposition}{0}
\setcounter{condition}{0}
\setcounter{assumption}{0}
\setcounter{conjecture}{0}
\setcounter{problem}{0}
\setcounter{remark}{0}
}

\newcommand{\noi}{\noindent}

\newcommand{\E}{\mathbb{E}}
\newcommand{\R}{\mathbb{R}}

\newcommand{\N}{\mathbb{N}}

\newcommand{\la}{\lambda}

\newcommand{\eps}{\varepsilon}

\newcommand{\al}{\alpha}

\newcommand{\om}{\omega}

\newcommand{\Gam}{\mathnormal{\Gamma}}

\newcommand{\Q}{{\mathbb Q}}

\newcommand{\PP}{{\mathbb P}}

\newcommand{\calA}{{\cal A}}

\newcommand{\calC}{{\cal C}}
\newcommand{\calD}{{\cal D}}

\newcommand{\calF}{{\cal F}}

\newcommand{\calJ}{{\cal J}}

\newcommand{\calQ}{{\cal Q}}\newcommand{\calO}{{\cal O}}

\newcommand{\calX}{{\cal X}}

\newcommand{\skp}{\vspace{\baselineskip}}

\newcommand\iy{\infty}

\newcommand{\clc}{{\cal C}}

\newcommand{\clo}{{\cal O}}

\newcommand{\cld}{{\cal D}}

\newcommand{\kaboom}{\hat\eps}

\title{Brownian control problems for 
	a multiclass M/M/1 queueing problem with model uncertainty\thanks{To appear in {\it Mathematics of Operations Research}}}

\author{Asaf Cohen\thanks{Department of Statistics,
		University of Haifa, 
		Haifa, 31905, Israel, 
shloshim@gmail.com
web: https://sites.google.com/site/asafcohentau/
}}

\date{\today}

\begin{document}

\maketitle

\begin{abstract} We consider a multidimentional Brownian control problem (BCP) with model uncertainty 
	that formally emerges from a multiclass M/M/1 queueing control problem under heavy-traffic with model uncertainty. 
	The BCP is formulated as a multidimensional stochastic differential game with two players: a minimizer that has an equivalent role to the decision maker in the queueing control problem and a maximizer whose role is to set up the uncertainty of the model. The dynamics are driven by a Brownian motion. We show that a state-space collapse propery holds. That is, the multidimensional 
	BCP can be reduced to a one-dimensional BCP with model uncertainty that also takes the form of a two-player stochastic differential game. Then, the value function of both games is characterized as the unique solution to a free-boundary problem from which we extract equilibria for both games. Finally, we analyze the dependence of the value function and the equilibria on the ambiguity parameters. 

\skp

\skp

\noi{\bf AMS Classification:} Primary: 93E20, 60K25, 91A15, 60J60; secondary: 49J15, 35R35.

\noi{\bf Keywords:} Brownian control problem; ambiguity aversion; model uncertainty; multiclass M/M/1; heavy-traffic; the Harrison--Taksar free-boundary problem. 
\end{abstract}

\section{Introduction}\label{sec_1}
\beginsec
Typically, heavy-traffic {\it queueing control problems} (QCPs) in the diffusion scale are treated by defining a limiting control problem associated with Brownian motion, called {\it Brownian control problems} (BCPs), first introduced by \cite{har1988}; for further reading on BCPs see e.g., \cite{bell-will-1,BG2006,BG2012} and the references therein. In this paper we study two BCPs with model uncertainty 
that formally emerge from a multiclass M/M/1 QCP with finitely many buffers with finite capacity under heavy-traffic with model uncertainty. The asymptotic relationship between the QCP and the BCPs is the subject of a separate paper, see \cite{Cohen2017}. Both the QCP and the BCP were studied in the classical case under the framework of G/G/1 without ambiguity about the model in \cite{ata-shi} and in \cite{PKH}, where the latter considers time constraints instead of the finite buffers constraints. These problems are referred to as {\it risk-neutral} problems.

\skp\noi
We are formulating a decision maker (DM) that  has a reference model in mind, which, up to some degree, describes the situation  she is facing. Since the DM is uncertain about the true model (either because the parameters are not fully known, can change over time, etc.), she takes into account other models and penalizes them based on their deviation from a specific reference model. The penalization depends then on how averse the DM is to ambiguity. Such ambiguity models are sometimes referred to as {\it model uncertainty} or {\it Knightian uncertainty}, see e.g., \cite{maenhout2004robust,Hansen2006,han-sar,bay-zha} and in the context of queueing systems see \cite{Shanti,blanchet2014robust,MR3544795}. 

\skp\noi
The first BCP with model uncertainty introduced is a {\it multidimensional stochastic differential game} (MSDG) with two players: a DM and the nature, which according to their goals are referred to as the {\it minimizer} and the {\it maximizer}, respectively. Borrowing terminology and 
the roles of the processes from the QCP to the MSDG (here and in the sequel), the state-space of the game is a product set of intervals of the form $[0,\hat b_i]$, where $\hat b_i$ is the capacity of buffer $i$. The minimizer in this game controls the server's effort allocation among the buffers and the admission/rejection to each buffer. The maximizer chooses the underlying probability measure; this is shown to be equivalent to stochastically perturbing the drift of the Brownian motion (possibly, differently for each coordinate).
The game's cost consists of holding and rejection penalties and a variant of the Kullback--Leibler divergence with respect to (w.r.t.) the relevant measures in this setup. The latter component stands for a penalty for the maximizer for changing the drift. 

\skp\noi
We show that a {\it state-space collapse} property holds. That is, we provide a one-dimensional BCP with uncertainty that also takes the form of a stochastic game, called the {\it reduced stochastic differential game} (RSDG), whose state emerges from the workload process in the QCP. The roles of the two players remain the same as in the MSDG and the dynamics and the cost functions have similar components. We show that the games are equivalent in the sense that given any strategy of the minimizer in either one of the games, we construct a strategy for the minimizer in the other game that performs at least as well, and therefore, also the value functions are the same (Proposition \ref{prop_31}). For further reading about workload reduction, see \cite{har1988,Har-Van,har2000,har2003,Har-Wil}. The advantage of such a reduction is that the dynamics live in a lower dimension and have only two components of singular controls that represent idleness and workload rejection. Therefore, most of the analysis is performed in the RSDG\footnote{Throughout the paper we alternatively refer to the BCPs as the `stochastic differential games', or as `the games', or explicitly by MSDG and RSDG.} setup. We characterize the value function of the RSDG as the unique classical solution of a {\it Hamilton-Jacobi-Bellman} (HJB) equation (Theorem \ref{thm_41}). Doing so, we extend the relationship between the reduced BCP and a relevant HJB equation studied in \cite[Equation (1.2)]{har-tak} and in \cite[Equation (41)]{ata-shi} to a similar relationship in a stochastic game setup with a different HJB; unlike in \cite{ata-shi}, due to the existence of a maximal player in our model, the HJB is not linear. Therefore, one cannot use the results given in these papers and rather needs to establish the relationship between the two. As a first step, we reduce the problem of solving the HJB equation to a free-boundary problem. (To the best of my knowledge, this is the first model uncertainty control problem that leads to a free-boundary problem with Neumann boundary conditions). Then we use the {\it shooting method} to solve the latter problem. In short, this is a method for solving boundary value problems using initial value problems. We take it one step forward in the free-boundary setup; see \cite[Section 7.3]{Stoer1980} for further reading about the method. 
%
Moreover, we supply equilibria in both games and refer to the equilibrium strategies of the minimizer (in both games) as the {\it  optimal strategy} of the minimizer (Theorem \ref{thm_42}). Starting with the RSDG, we show that the optimal strategy for the minimizer in this game is a reflecting strategy. Namely, the minimizer should use minimal idleness and minimal amount of workload rejections in order to keep the workload in a fixed interval. The equilibrium strategy of the maximizer in this game is also provided. From the RSDG equilibrium we construct an equilibrium in the MSDG. According to the minimizer's optimal policy in the MSDG, the queue length processes evolve along a certain curve in the state-space. The minimizers' opimal strategies in our games are shown to have the same structures as the optimal strategies in the risk-neutral BCPs in \cite{ata-shi}. Such a result is not so obvious for reasons having to do with the non-stationarity structure of the problem caused by the existence of the maximizer player. 
The difference between the policy of the minimizer in the multidimensional setup here as compared to \cite{ata-shi} 
is only the cut-off level of the reflecting strategy in the RSDG, which affects the point of reflection on the multidimensional curve. 

\skp\noi
Aside from studying the games, we also analyze the dependence of the games on the ambiguity parameters. We show continuity of the value function and the optimal reflecting strategies w.r.t.~the ambiguity parameters and that as the ambiguity vanishes, the problem converges to the risk-neutral problem studied in \cite{ata-shi} (Theorems \ref{thm_45} and \ref{prop_46}).

\skp\noi
We now discuss about the position of the current work between the risk-neutral BCPs and the deterministic differential games (DDGs) studied in \cite{ata-coh} and provide future outlook. 
Recall that the risk-neutral BCPs from \cite{ata-shi} governs the limiting behavior of a diffusion scaled multiclass G/G/1 QCP with linear utility function. On the other side, the DDGs from \cite{ata-coh} approximate the same type of QCP with the following differences: the moderate-deviation replaces the diffusion scaling and the utility function is exponential, scaled with the moderate-deviation parameters. While in the risk-neutral case the probability of buffer overflow is approximated by the probability that a Brownian motion with drift hits a positive level, which is of order $\clo(1)$, in the moderate-deviation scaling this probability vanishes with the scaling parameter. As a compensation, the criteria considered is a risk-sensitive one, which gives an overwhelming scaled exponential cost for such event. This means that the DM has a large ambiguity about the model and  she gives high weight even for events that are `very' unlikely to happen. The intuition behind the relationship between the QCP in the moderate-deviation heavy-traffic regime and the limiting DDG goes back to the classical risk-sensitive control problem with small noise diffusion. As argued in \cite[Ch.~VI.2]{Fleming2006}, consider a cost function $\calJ$ and a positive constant $\varrho$, then the risk-sensitive cost can be expressed as,
\begin{align}\label{newasaf1049}
	\al^{-1}\log\left(E^\PP[e^{\al\calJ}]\right)=\sup_{\Q}E^\Q\Big[\al\Big(\calJ-\log\left(\frac{d\Q}{d\PP}\right)\Big)\Big],
\end{align}
where the supremum is taken over an appropriate set of measures. When the noise coefficient is $\al^{-1/2}$ and $\al\to \iy$, the limiting problem is a DDG. What we consider is the criteria given in \eqref{newasaf1049}, without sending $\al\to \iy$. In the QCP, it means that we still consider the classical diffusion scaling, yet our criteria gives space for ambiguity about the true underlying probability measure. Since the higher $\al$ is, the more weight the DM gives to unlikely events, we say that the DM has more ambiguity in the moderate-deviation risk sensitive QCP then in the one related to our work, and thus we position our model between the two models described above. This work does not aim to establish a deeper connection between the model uncertainty's BCPs and the DDGs. At this point it worth mentioning that the relation between the Kullback--Leibler-constrained formulation and the exponential functional (i.e., \eqref{newasaf1049}) is well known in other contexts, see e.g., \cite{chowdhary2013distinguishing} for an example in the context of uncertainty quantification.

\skp\noi
In \cite{ata-coh} it was shown that the optimal strategies of the minimizers in the DDGs have the same structure as the ones in the BCPs from \cite{ata-shi}. Together with the observation stated earlier, we get that all the three models (\cite{ata-shi}, \cite{ata-coh}, and the present one) share the same optimal policy, where the different lies in the cut-off point for rejections. Moreover, in a recent line of research of queueing systems under the moderate-deviation heavy-traffic regime with risk-sensitive performance criteria it was shown that classical results from the theory of risk-neutral QCPs and BCPs such as state-space collapse and generalyzed $c\mu$ rule hold in these models as well, see e.g., \cite{ata-bis,bis2014,ata-coh,ata-sub_cmu,ata-coh2017, ata-men}. The current paper together with \cite{Cohen2017} represents a line of research of QCPs and their associated BCPs under model uncertainty. 

\skp\noi
In this work we assume that the buffers are finite. This property arises naturally in private cloud computing, which are limited in data and capacity. For a more detailed application, see \cite{shi-ata-cid}. Intuitively, the problem with unlimited capacities and without rejections seems to be simpler. We expect that similar tools can be used to show that the optimal strategy would be to assign {\it fixed} lowest priority to the class with the smallest $\hat h_i\mu_i$ value. Consequently, the limiting queueing sizes would vanish with the scaling parameter for the rest of the buffers. This is left for future work.

\skp\noi
In summary, our main contributions are as follows. We
\begin{itemize} \itemsep0em
	\item provide and solve for the first time a BCP with model uncertainty that emerges from 
	a multiclass M/M/1 QCP with model uncertainty;
	\item show that a state-space collapse property holds for this game (Corollary \ref{cor_31});
	\item show that the reduced game solves uniquely a relevant HJB equation, which is a nonlinear free-boundary problem and that there is an optimal reflecting strategy (Theorem \ref{thm_41});
	\item provide equilibria for the two games considered (Theorem \ref{thm_42});
	\item analyze the dependence of the value function and the equilibria on the ambiguity parameters (Theorems \ref{thm_45} and \ref{prop_46}).
\end{itemize}

\skp\noi
The paper is organized as follows. In Section \ref{sec_2} 
we motivate and present the stochastic differential games and study the relationship between the two. Next, In Section \ref{sec_4} we study the RSDG. We provide the HJB equation and prove that the value function of the game is the unique smooth solution of the HJB. Moreover, we show that the minimizer has an optimal reflecting strategy. In Section \ref{sec_5} we discuss the uniqueness of the optimal reflecting strategy and find equilibria in both games. Finally, in Section \ref{sec_6} we study the dependency on the ambiguity parameters.

\subsection{Notation} 
We use the following notation.
For $a,b\in\R$, $a\wedge b=\min\{a,b\}$ and $a\vee b=\max\{a,b\}$. For a positive integer $k$ and $c,d\in\R^k$, $c\cdot d$ denotes the usual
scalar product and $\|c\|=(c\cdot c)^{1/2}$. We denote $[0, \iy)$ by $\R_+$. For subintervals $I_1,I_2\subseteq\R$ and $m\in\{1,2\}$ we denote by $\clc(I_1,I_2)$, $\clc^m(I_1,I_2)$, and $\cld(I_1, I_2)$ the space of continuous functions [resp., functions with continuous derivatives of order $m$, functions that are right-continuous with finite left limits (RCLL)] mapping $I_1\to I_2$. 
The space $\calD(I_1,I_2)$ is endowed with the usual Skorohod topology.

\section{The BCPs}\label{sec_2}
\beginsec

We start this section with a motivative QCP. The two BCPs are presented in Section \ref{sec_3}:                   in Section \ref{sec_31} we formally derive the MSDG and in Section \ref{sec_32} its reduction. Then in Section \ref{sec_33} we show that the two games share the same value and that given any strategy (for the minimizer) in either one of the games, one can construct a strategy in the second game that performs at least as well. 

\subsection{Motivative QCP}
Consider a model that consists of $I$ customer classes and a single server. Each class has its own finite buffer and upon arrival, customers are queued in the corresponding buffer or rejected. Within each class, customers are served at the order of their arrivals. Processor sharing is allowed and the server may serve up to $I$ customers at a time, where two customers from the same class cannot be served simultaneously. 
The system under consideation is in heavy-traffic. For this, we consider a sequence of systems, indexed by the scaling  parameter $n\in\N$. For every 
$n$ we consider a {\it reference probability space} that supports 
independent Poisson processes $A^n_i$ and $S^n_i$, $i\in[I]:=\{1,\ldots,I\}$ with rates $\lambda^n_i$ and $\mu^n_i$, repsectively. The value $A^n_i(t)$ stands for the number of customers of class $i$ that arrived to the system until time $t\in\R_+$, and  
$S^n_i(t)$ is the number of service completions of class $i$ customers had the server dedicated all of its effort to class $i$ during the time interval $[0,t]$.

Denote by $T^n_i(t)$ the units of time that the server devoted to class $i$ until time $t$. 
For every $t\in\R_+$ and $i\in[I]$, $S^n_i(T^n_i(t))$ is the number of service completions of class $i$ customers until time $t$. This is a Cox process with the infinitesimal intensity $\mu^n_idT^n_i(t)$. 
Rejections of customers are allowed upon arrival and a rejected customer will never return to the system. The number of customers from class $i$ that were rejected by time $t$ is denoted by $R^n_i(t)$. For every $i\in[I]$, the balance equation is given by,
\begin{equation}\label{204a}
	X^n_i(t)=X^n_i(0)+A^n_i(t)-S^n_i(T^n_i(t))-R^n_i(t),\qquad t\in\R_+,
\end{equation}
where $X^n_i(t)$ stands for the number of class $i$ customers in the system at time $t$. 
We use the notation $L^n=(L^n_i)_{i=1}^I$ for $\{A,S,X,R,T\}$.

We assume that
\begin{align}\label{205}
	\lambda^n_i:=\lambda_i n+\hat\lambda_in^{1/2}+o(n^{1/2}),\qquad
	\mu^n_i:=\mu_in+\hat\mu_in^{1/2}+o(n^{1/2}),
\end{align}
where $\lambda_i,\mu_i\in(0,\iy)$ and $\hat\lambda_i,\hat\mu_i\in\R$ are fixed. Moreover, the system is assumed to be \emph{critically loaded}, that is, 
$\sum_{i=1}^I\rho_i=1$, where $\rho_i:=\lambda_i/\mu_i$, $i\in[I]$.

The scaled version of \eqref{204a} is given by,
\begin{equation} \label{newA}
	\hat{X}^n_i(t)=\hat{X}^n_i(0) + \hat m^n_it
	+\hat{A}^n_i(t)-\hat{S}^n_i(T^n_i(t))
	+\hat  Y^n_i(t)- \hat{R}^n_i(t),\qquad t\in\R_+,
\end{equation}
where
\begin{align}\label{207}
	&  \hat X^n_i(t):=n^{-1/2}X^n_i(t),\quad \hat{A}^n_i(t):=n^{-1/2}(A^n_i(t)-\la^n_i t),\quad
	\hat{S}^n_i(t):=n^{-1/2}(S^n_i(t)-\mu^n_i t),\\\notag
	&\quad\qquad\hat Y^n_i(t):=\mu^n_in^{-1/2}(\rho_i t-T^n_i(t)), \qquad \hat R^n(t):=n^{-1/2}R^n(t),
\end{align}
and
\begin{equation}\notag
	\hat m^n_i:=n^{-1/2}(\lambda^n_i-\rho_i\mu^n_i).
\end{equation}
As previously, we use the notation 
$\hat L^n=(\hat L^n_i)_{i=1}^I$  for $L\in\{A,S,X,R,Y,T,m\}$.

The capacity of buffer $i$ is given by $ \hat b^n_i:=\hat b_in^{1/2}$ for some constant $\hat b_i\in(0,\iy)$, $i\in[I]$. We assume that $\hat X^n(0)\in\calX:=\prod_{i=1}^I[0, \hat b_i]$, and the rejection mechanism assures that 
\begin{align}\label{208}
	\hat X^n_i(t)\in\calX,\qquad t\in\R_+,\;\; \PP^n\text{-a.s.}
\end{align}
We now present the optimization criteria. The intuition behind it is as follows. The DM, also referred to as the {\it minimizer}, chooses a control based on the past observations. He minimizes a cost that takes into account a possible deviation from the reference model. For this, we consider an adverse `player', also referred to as the {\it maximizer}, who has access to the policy chosen by the minimizer and to the history. 
This player is penalized for deviating from the reference model.

In details, the vectors $\hat h,\hat r\in(0,\iy)^I$ stand for the holding and rejection costs, respectively. 
%
The DM is uncertain about the underlying reference probability measure, or in other words,  she suspects with some level of uncertainty that the rates of the processes $\{A^n_i\}_{i=1}^I$ and $\{S^n_i\}_{i=1}^I$ are not exactly $\{\la^n_i\}_{i=1}^I$ and $\{\mu^n_i\}_{i=1}^I$ and may be unspecified or may even change over the time. Therefore, she considers a set of candidate measures and penalizes their deviation from the reference probability measure. The penalization is done by using discounted variants of the Kullback--Leibler divergence, given by
\begin{align}\label{211}
	\begin{split}
		L^\varrho_1(\hat \Q^n_{1,i}\|\PP^n_{1,i})&:=
		\E^{\hat \Q^n_{1,i}}\left[\int_0^\iy\varrho e^{-\varrho t}\log\frac{d\hat \Q^n_{1,i}(t)}{d\PP^n_{1,i}(t)}dt\right], \qquad n\in\N,\;i\in[I],\\
		L^\varrho_2(\hat \Q^n_{2,i}\|\PP^n_{2,i})&:=
		\E^{\hat \Q^n_{2,i}}\left[\int_0^\iy\varrho e^{-\varrho t}\log\frac{d\hat \Q^n_{2,i}(t)}{d\PP^n_{2,i}(t)}dT^n_i(t)\right], \qquad n\in\N,\;i\in[I],
	\end{split}
\end{align}
where $\PP^n_{1,i}$ and $\PP^n_{2,i}$ are the reference measures under which $A^n_i$ and $S^n_i$ are Poisson processes with rates $\la^n_i$ and $\mu^n_i$, respectively. 
To establish the level of ambiguity, for every $i\in[I]$, we consider the (finite and) positive parameters $\kappa_{1,i}$ and $\kappa_{2,i}$ that quantify the amount of ambiguity that the DM has regarding the rates $\la^n_i$ and $\mu^n_i$, or in other words, the measures $\PP^n_{1,i}$ and $\PP^n_{2,i}$, respectively. 
Set $\kappa:=(\kappa_{1,i},\kappa_{2,i})_{i=1}^I$. The DM is facing the following optimization problem:
\begin{align}\notag
	V^n( X^n(0);\kappa)= \;&
	\inf_{(T^n,R^n)}\;\sup_{\hat{\Q}^n\in\hat{ \cal{Q}}^n(\hat X^n(0))} \hat J^n(\hat X^n(0), U^n, R^n,\hat\Q^n;\kappa),
\end{align}
where
\begin{align}\label{212}
	& \hat J^n( \hat X^n(0),U^n,R^n,\hat\Q^n;\kappa)=\\\notag
	&\quad\E^{\hat \Q^n}\Big[\int_0^\iy e^{-\varrho t}\left(\hat h\cdot\hat  X^n(t)dt +\hat r\cdot  d\hat R^n(t)\right) \Big]
	-\sum_{i=1}^I\frac{1}{\kappa_{1,i}}L^\varrho_1(\hat \Q^n_{1,i}\|\PP^n_{1,i})
	-\sum_{i=1}^I\frac{1}{\kappa_{2,i}}L^\varrho_2(\hat \Q^n_{2,i}\|\PP^n_{2,i}),
\end{align}
$\hat\Q^n=\prod_{i=1}^I(\hat\Q^n_{1,i}\times\hat\Q^n_{2,i})$, and the set of candidate measures $\hat{\cal{Q}}^n(\hat X^n(0))$ is described at the end of this paragraph. 
When $\kappa_{j,i}$ is `small' (resp., `large') we say that there is a weak (resp., strong) ambiguity about the rates of the processes $A^n_i$ and $S^n_i(T^n_i):=S^n_i(T^n_i(\cdot))$. 
The idea is that for small $\kappa_{j,i}$'s there is a big punishment per unit of deviation from the reference measure and therefore, the measures
$\hat\Q^n_{j,i}$ and $\PP^n_{j,i}$ should be close to each other and as a consequence also the relevant expectations. 
We now turn to define the set of candidate measures. 
A probability measure $\hat\Q^n$ belongs to $\hat{\cal{Q}}^n(\hat X^n(0))$ if for every $i\in[I]$ and $ t\in\R_+$ it
satisfies
\begin{align}\label{216}
	\frac{d\hat \Q^n_{1,i}(t)}{d\PP^n_{1,i}(t)}&=\exp\Big\{\int_0^t\log\left(\frac{\psi^n_{1,i}(s)}{\la^n_i}\right) d A^n_i(s)-\int_0^t (\psi^n_{1,i}(s)-\la^n_i)ds\Big\}, \\\label{216z}
	\frac{d\hat \Q^n_{2,i}(t)}{d\PP^n_{2,i}(t)}&=\exp\Big\{\int_0^t\log\left(\frac{\psi^n_{2,i}(s)}{\mu^n_i}\right) d S^n_i(T^n_i(s))-\int_0^t (\psi^n_{2,i}(s)-\mu^n_i)dT^n_i(s)\Big\},
\end{align}
for measurable and positive processes $\psi^n_{j,i}$ that are predictable w.r.t.~the filtration generated by the arrival and service completions processes, satisfying $\int_0^t\psi^n_{j,i}(s)ds<\iy$ $\PP^n_{j,i}$-a.s., for every $t\in\R_+$. These conditions guarantee that the right-hand sides in \eqref{216} are $\PP^n_{j,i}$-martingales, and that under the measure $\hat\Q^n_{1,i}$ (resp., $\hat\Q^n_{2,i}$), the processes $A^n_i$ (resp., $S^n_i(T^n_i)$) is a counting process with infinitesimal intensity $\psi^n_{1,i}(t)dt$ (resp., $\psi^n_{2,i}dT^n_i(t)$). Notice that we do not force the critically loaded condition under the measures $\hat\Q^n_{j,i}$. As we argue in \cite[Section 4]{Cohen2017} such changes of measures are `too costly' and will be avoided by the maximizer who would choose 
$\psi^n_{1,i}(t)=\la_i^n+\calO(n^{1/2})$ and $
\psi^n_{2,i}(t)=\mu_i^n+\calO(n^{1/2}).$

We now provide an approximation to the change of measure penalty that will be useful as a motivation in the next section. 
Since $A^n_i(\cdot)-\int_0^\cdot\psi^n_{1,i}(s)ds$ is a martingale, we have
\begin{align}\notag
	\begin{split}
		&L^\varrho_1(\hat \Q^n_{1,i}\|\PP^n_{1,i})\\
		&\quad=\E^{\hat\Q^n_{1,i}}\Big[\int_0^\iy\rho e^{-\varrho t}\Big(\int_0^t\log\Big(\frac{\psi^n_{1,i}(s)}{\la^n_i}\Big)dA^n_i(s)-\int_0^t(\psi^n_{1,i}(s)-\la^n_i)ds\Big)dt\Big]\\
		&\quad=\E^{\hat\Q^n_{1,i}}\Big[\int_0^\iy\rho e^{-\varrho t}\Big(\int_0^t\log\Big(\frac{\psi^n_{1,i}(s)}{\la^n_i}\Big)(dA^n_i(s)-\psi^n_{1,i}(s)ds)\\
		&\qquad\qquad\qquad\qquad\qquad\qquad+\int_0^t\Big\{\psi^n_{1,i}(s)\log\Big(\frac{\psi^n_{1,i}(s)}{\la^n_i}\Big)-\psi^n_{1,i}(s)+\la^n_i\Big\}ds\Big)dt\Big]\\
		&\quad=\E^{\hat\Q^n_{1,i}}\Big[\int_0^\iy\rho e^{-\varrho t}\Big(\int_0^t\Big\{\psi^n_{1,i}(s)\log\Big(\frac{\psi^n_{1,i}(s)}{\la^n_i}\Big)-\psi^n_{1,i}(s)+\la^n_i\Big\}ds\Big)dt\Big]\\
		&\quad=\E^{\hat\Q^n_{1,i}}\Big[\int_0^\iy\rho e^{-\varrho t}\Big\{\psi^n_{1,i}(t)\log\Big(\frac{\psi^n_{1,i}(t)}{\la^n_i}\Big)-\psi^n_{1,i}(t)+\la^n_i\Big\}dt\Big]
		,
	\end{split}
\end{align}
where a change of variables is used to obtain the last equality. The same analysis applies for $S^n_i$ as well. 
Consider the processes
\begin{align}\notag
	\hat \psi^n_{1,i}(t):=( \la_i n)^{-1/2}(\psi^n_{1,i}(t)-\la^n_i),\qquad
	\hat \psi^n_{2,i}(t):=(\mu_i n)^{-1/2}(\psi^n_{2,i}(t)-\mu^n_i).
\end{align}
The Taylor expansion of $\log (1+x)$, gives\footnote{This is rigorously justified in \cite{Cohen2017} by another level of approximation to the processes $\{\psi^n_{j,i}\}$.}
\begin{align}\label{ac08}
	&\frac{1}{\kappa_{1,i}}L^\varrho_1(\hat \Q^n_{1,i}\|\PP^n_{1,i})+\frac{1}{\kappa_{2,i}}L^\varrho_2(\hat \Q^n_{2,i}\|\PP^n_{2,i})\\\notag
	&\qquad\approx \E^{\hat\Q^n_{1,i}\times\hat\Q^n_{2,i}}\left[\int_0^\iy e^{-\varrho t}\left\{\frac{1}{2\kappa_{1,i}} (\hat\psi_{1,i}^n(t))^2+\frac{1}{2\kappa_{2,i}}\rho_i (\hat\psi_{2,i}^n(t))^2\right\}dt\right].
\end{align}
The term $\rho_i$ is due to the convergence $T^n_i(t)\to\rho_it$, $t\in\R_+$.  
Now, since the maximizer is free to choose $\hat\psi^n_{1,i}$ and $\hat\psi^n_{2,i}$,  she faces the two steps optimization problem. First, to choose $\hat\psi^n_i(t)$ and then to solve
\begin{align}\notag
	\min_{(\hat\psi_{1,i}^n(t),\hat\psi_{2,i}^n(t))}\;\left\{\frac{1}{2\kappa_{1,i}} (\hat\psi_{1,i}^n(t))^2+\frac{1}{2\kappa_{2,i}}\rho_i (\hat\psi_{2,i}^n(t))^2:\la_i^{1/2}  \hat\psi^n_{1,i}(t)-\rho_i\mu_i^{1/2}  \hat\psi^n_{2,i}(t)=(2\la_i)^{1/2}\hat\psi^n_i(t)\right\}.
\end{align}
The minimal value of the above equals $(2\kaboom_i)^{-1} (\hat\psi_i^n(t))^2$, where
\begin{align}\label{300b}
	\kaboom_i:=\frac{1}{2}(\kappa_{1,i}+\kappa_{2,i}).
\end{align}
Therefore, 
\eqref{ac08} can be approximated by $(2\kaboom_i)^{-1} (\hat\psi_i^n(t))^2$. 
These arguments are rigorously justified in \cite[Section 4]{Cohen2017}.

\subsection{Two stochastic differential games}\label{sec_3}

\subsubsection{The multidimensional stochastic differential game (MSDG)}\label{sec_31}

The sequence of the scaled and centered $2I$-dimensional Poisson processes $(\hat A^n,\hat S^n)$ weakly converges to a $2I$-dimensional Brownian motion starting at zero, with zero mean and the covariance matrix 
$
\text{Diag}\;(\la_1^{1/2} ,\ldots,\la_I^{1/2},\mu_1^{1/2} ,\ldots,\mu_I^{1/2}). 
$ Formally speaking, if the process $\hat Y^n$ is of order one as $n\to\iy$, which is rigorously proven in \cite[Section 4]{Cohen2017}, we get from its definition in \eqref{207} that $T^n(t)\to(\rho_1t,\ldots,\rho_It)$, $t\in\R$, and therefore, under $\PP^n$, $(\hat A^n_i-\hat S^n_i(T^n_i))_{i=1}^I$ weakly converges to an $I$-dimensional Brownian motion starting at zero, with zero mean and the covariance matrix
\begin{align}\notag
	\hat\sigma=(\hat\sigma_{ij}):=\text{Diag}\left((2\la_1)^{1/2} ,\ldots,(2\la_I)^{1/2}\right),
\end{align}
where $\hat S^n_i(T^n_i):=\hat S^n_i(T^n_i(\cdot))$, $i\in[I]$.

Recall that in the QCP, an admissible control is of the form $(U^n,R^n)$. Notice that $(\hat Y^n(t),\hat R^n(t))$ is uniquely determined by $(U^n(s),R^n(s))_{0\le s\le t}$. In the MSDG we consider two $I$-dimensional processes, $\hat R$ and $\hat Y$ that play the roles of instantaneous controls, which stand for the scaled rejection process $\hat R^n$ and the scaled idle time process $\hat Y^n$, respectively. 
Moreover, the ambiguity about the true underlying probability measure in the QCP, which is formulated by a penalty for deviating from the reference measure is translated to the limiting problem as well. In the definition below we refer to two players by their roles as a minimizer and a maximizer even though the roles can only be derived from the cost function, which is presented afterwards. 

\begin{definition}[admissible controls, MSDG]\label{def_32} 
	An {\it admissible control for the minimizer} for any initial state $\hat x_0\in\calX$ is a filtered probability space  
	\begin{align}\notag
		(\Omega,\calF,\{\calF_t\},\PP):=\Big(\prod_{i=1}^I\Omega^i,\calF^1\otimes\ldots\otimes\calF^I,\{\calF_t\},
		\prod_{i=1}^I\PP_i\Big),
	\end{align}
	that supports a process $(\hat Y,\hat R)$ taking values in $(\R^I_+)^2$ with RCLL sample paths adapted to the filtration $\{\calF_t\}$, where $(\Omega^i,\calF^i,\{\calF^i_t\},\PP_i)$ supports a one-dimensional standard Brownian motion $\hat B_i$ adapted to the filtration $\{\calF^i_t\}$, $i\in[I]$. Moreover, assume that the following properties hold:
	

	\noi
	(i)
	\begin{align}\label{300}
		\text{for every $i\in[I]$ and $0\le s<t$, $\hat B_i(t)-\hat B_i(s)$ is independent of $\calF^i_s$ under $\PP_i$;}
	\end{align}
	
	\noi(ii)
	\begin{align}\label{301}
		&\text{ $\theta\cdot\hat  Y$ and $\hat R_i$, $i\in[I]$ are nonnegative and nondecreasing, where $\theta:=(\mu_1^{-1},\ldots,\mu_I^{-1})$};
	\end{align}
	
	\noi
	(iii) 
	\begin{align}\label{302}
		\hat X(t)=\hat x_0+\hat mt+\hat \sigma \hat B(t)+\hat Y(t)-\hat R(t),\quad t\in\R_+,
	\end{align}
	such that 
	\begin{align}\label{303}
		\hat X(t)\in\calX,\quad t\in\R_+,\;\PP\text{-a.s.,}
	\end{align}
	where $\hat m:=\lim_n\hat m^n=(\hat\la_i-\rho_i\hat\mu_i)_{i=1}^I$ and $\hat B=(\hat B_i)_{i=1}^I$.
	
	An {\it admissible control for the maximizer} is a product measure $\hat \Q=\prod_{i=1}^I \hat\Q_i$, where each $\hat\Q_i$ is defined on $(\Omega^i,\calF^i,\{\calF^i_t\})$, such that 
	\begin{align}\label{304}
		\frac{d\hat \Q_i(t)}{d\PP_i(t)}=\exp\Big\{\int_0^t\hat \psi_i(s) d\hat B_i(s)-\frac{1}{2}\int_0^t \hat \psi^2_i(s)ds\Big\},\quad t\in\R_+,
	\end{align}
	for an $\{\calF_t\}$-progressively measurable process $\hat \psi=(\hat\psi_1,\ldots,\hat\psi_I)$ satisfying 
	\begin{align}\label{305}
		&\E^{\PP}\Big[\int_0^\iy e^{-\varrho s}\hat \psi_i^2(s)ds\Big]<\iy\quad\text{and}\quad\E^{\PP}\Big[e^{\frac{1}{2}\int_0^t\hat \psi^2_i(s)ds}\Big]<\iy\quad t\in\R_+,\;i\in[I].
	\end{align}
\end{definition}
We consider a probability space that is constructed from $I$ small probability spaces, where each one supports the processes associated with one of the classes (in the QCP studied in \cite{Cohen2017} there are $2I$ small probability spaces as can be inferred from the structure of $\hat\Q^n$). The Brownian motion approximates the difference $(\hat A^n_i-\hat S^n_i(T^n_i))_{i=1}^I$ up to the deterministic covariance matrix, and \eqref{302} follows by  \eqref{204a}. The emphasis of Condition \eqref{300} is that the independent condition for $\hat B_i$ is given w.r.t.~the filtration $\calF^i_t$ and not simply w.r.t.~its own filtration, which merely follows by being a Brownian motion. It is necessary for the approximation procedure, see the details in \cite[Section 4]{Cohen2017}. Condition \eqref{301} follows since the rejection process (in the QCP) and also $\theta^n\cdot\hat Y^n$ are nondecreasing. Occasionally, we refer to $\hat R$ as the rejection process in the MSDG. The buffer constraint is imposed in \eqref{303}. Pay attention that we consolidate the processes $\hat A^n_i$ and $\hat S^n_i(T^n_i)$ into one Brownian motion. Hence, we consider only $I$ changes of measures instead of $2I$. 

We now explain the change of measure structure. Recall the definition of $\hat \psi^n_{j,i}$ and $\hat \psi_i$ from the previous section. 
Then,
\begin{align}\label{218e}
	&\hat{A}^n_i (t)
	=n^{-1/2}\Big(A^n_i(t)-\int_0^t\psi^n_{1,i}(s)ds\Big)
	+\la_i^{1/2}\int_0^t\hat\psi^n_{1,i}(s)ds,\\\notag
	&\hat{S}^n_i (T^n_i(t))
	=n^{-1/2}\Big({S}^n_i (T^n_i(t))-\int_0^t\psi^n_{2,i}(s)dT^n_i(s)\Big)
	+\mu_i^{1/2}\int_0^t\hat\psi^n_{2,i}(s)dT^n_i(s).
\end{align}
Informally speaking, in both lines above, under $\hat \Q^n_{1,i}\times\hat \Q^n_{2,i}$, the first term is approximately a standard Brownian motion and the second term approximates a drift (in case it converges). This claim is handelled rigorously in \cite[Section 4]{Cohen2017}. As a result, under $\hat\Q^n_i$, $\hat A^n_i-\hat S^n_i(T^n_i)$ is approximately a diffusion process with drift $(2\la_i)^{1/2}\hat\psi^n_i:= \la_i^{1/2}  \hat\psi^n_{1,i}-\rho_i\mu_i^{1/2}  \hat\psi^n_{2,i}$ and a diffusion coefficient $\sigma_{ii}=(2\la_i)^{1/2}$; the term $\rho_i$ is due to the convergence $T^n_i(t)\to\rho_it$, $t\in\R_+$.   

%
%

\begin{remark}\label{rem_31}
	(i) Notice that the structure of the information in the game is consistent with the one in the QCP. The minimizer chooses a strategy and the maximizer, which is penalized for deviating from the reference measure, responds to this strategy by choosing a worst case scenario. For further reading about the structure of the information in control problems with model uncertainty, the reader is referred to \cite{Sirbu2014}.
	
	(ii) Given any $\{\calF_t\}$-progressively measurable process $\hat\psi$ that satisfies the conditions in \eqref{305}, the right-hand side (r.h.s.) of \eqref{304} is a martingale, and therefore, there exists a probability measure $\hat\Q$ such that $\hat\Q|_{\calF_t}$ satisfies \eqref{304} for all $t\in\R_+$. 
	
	(iii) Equation \eqref{302} can alternatively be written as
	\begin{align}\label{306}
		\hat X(t)=\hat x_0+\hat mt+\int_0^t\hat\sigma\hat\psi(s)ds+\hat \sigma \hat B^{\hat\Q}(t)+\hat Y(t)-\hat R(t),\quad t\in\R_+,
	\end{align}
	where $\hat B^{\hat \Q}(t):=\hat B(t)-\int_0^t\hat\psi(s)ds$, $t\in\R_+$, is an $\{\calF_t\}$-$I$-dimensional standard Brownian motion under $\hat\Q$.
	
\end{remark}

Denote by $\hat \calA(\hat x_0)$ the set of all admissible controls for the minimizer, given the initial condition $\hat x_0$. We often abuse notation and denote $(\hat Y,\hat R)\in\hat \calA(\hat x_0)$, keeping in mind that the control includes a filtered probability space. The set of all admissible controls for the maximizer is denoted by $\hat \calQ(\hat x_0)$. 

\skp\noi{\bf The cost function (MSDG).} 
Recall the discussion about the approximation of the cost function of the QCP from the previous section and \eqref{300b}. 
Set $\kaboom=(\kaboom_i)_{i=1}^I$. The cost associated with the initial condition $\hat x_0$ and the strategies $(\hat Y,\hat R)$ and $\hat \Q$ is given by
\begin{align}\label{307}
	\hat J(\hat x_0,\hat Y,\hat R,\hat \Q;\kaboom):= \;&
	\E^{\hat \Q}\Big[\int_0^\iy e^{-\varrho t}\left(\hat h\cdot\hat  X(t)dt +\hat r\cdot  d\hat R(t)\right) \Big]-\sum_{i=1}^I\frac{1}{\kaboom_i}L^\varrho(\hat \Q_i\|\PP_i),
\end{align}
where 
\begin{align}\label{307a}
	L^\varrho(\hat \Q_i\|\PP_i):=
	\E^{\hat \Q_i}\left[\int_0^\iy\varrho e^{-\varrho t}\log\frac{d\hat \Q_i(t)}{d\PP_i(t)}dt\right] 
\end{align}
The cost function can alternatively be expressed by
\begin{align}\label{308}
	\hat J(\hat x_0,\hat Y,\hat R,\hat \Q;\kaboom)=\;&
	\E^{\hat \Q}\Big[\int_0^\iy e^{-\varrho t}\Big(\hat h\cdot\hat  X(t)dt +\hat r\cdot  d\hat R(t)-\sum_{i=1}^I\frac{1}{2\kaboom_i}\hat \psi^2_i(t)dt\Big) \Big],
\end{align}
with $\hat \psi$ satisfying \eqref{304}--\eqref{305} above. Indeed, 
\begin{align}\label{309}
	L^\varrho(\hat \Q_i\|\PP_i)&=\E^{\hat\Q_i}\left[\int_0^\iy \varrho e^{-\varrho t} \left(-\frac{1}{2}\int_0^{t}\hat\psi^2_i(s)ds+\int_0^{t}\hat\psi_i(s) d\hat B_i(s)\right)dt\right]\\\notag
	&=\E^{\hat\Q_i}\left[\int_0^\iy \varrho e^{-\varrho t} \left(-\frac{1}{2}\int_0^{t}|\hat\psi_i(s)|^2ds+\int_0^{t}\hat\psi_i(s)\cdot(d\hat B_i^{\hat\Q}(s)+\hat\psi_i(s)ds)\right)dt\right]\\\notag
	&=\E^{\hat\Q_i}\left[\int_0^\iy \varrho e^{-\varrho t} \left(\frac{1}{2}\int_0^{t}\hat\psi_i^2(s)ds\right)dt\right]
	=\E^{\hat\Q_i}\left[\frac{1}{2}\int_0^\iy\hat\psi_i^2(s)\int_s^\iy \varrho e^{-\varrho t}dtds\right]\\\notag
	&=\E^{\hat\Q_i}\left[\frac{1}{2}\int_0^\iy e^{-\varrho t}\hat\psi_i^2(t)dt\right]<\iy.
\end{align}
Compare this structure with the approximated penalty in the QCP given at the end of the previous section.

While the form of the cost function given in \eqref{307} captures better the ambiguity aversion, the form of the cost given in \eqref{308} is more useful from a technical point of view. Moreover, the dynamics in \eqref{306} together with the cost function given in \eqref{308} are similar in their structure to their correspondences in Equation (11) and the display below (13) together with (2) in \cite{ata-coh}. The DM is faced the following optimization problem
\begin{align}\notag
	\hat V(\hat x_0;\kaboom)=\inf_{(\hat Y,\hat R)\in\hat \calA(\hat x_0)}\;\sup_{\hat \Q\in\hat \calQ(\hat x_0)}\;\hat J(\hat x_0,\hat Y,\hat R,\hat \Q;\kaboom)
\end{align}

\subsubsection{The reduced stochastic differential game (RSDG)}\label{sec_32}
We now present the RSDG. This game is one-dimensional and obtained by projecting the processes from \eqref{302} in the $\theta$ direction, which is given in \eqref{301}. For this we need the following notation, 
\begin{align}\label{310b}
	x_0:=\theta\cdot \hat x_0,\quad m:=\theta\cdot\hat  m,\quad\sigma:=\|\theta\hat\sigma\|,
\end{align}
\begin{align}\label{310g}
	\eps:=\frac{1}{\sigma^2}\sum_{i=1}^I(\theta\hat\sigma)_i^2\kaboom_i,
\end{align}
and
\begin{align}\notag
	b:=\max\{\theta\cdot \hat \xi:\hat \xi\in\calX\}=\theta\cdot\hat b,
\end{align}
where $\hat b=(\hat b_i)_{i=1}^I$.
\begin{definition} [admissible controls, RSDG]\label{def_33}
	An {\it admissible control for the minimizer} for any initial state $x_0\in[0,b]$ is a filtered probability space  $(\Omega,\calF,\{\calF_t\},\PP)$ that supports a one-dimensional standard Brownian motion $ B$ and a process $ ( Y, R)$ taking values in $\R_+^2$ with RCLL sample paths, both adapted to the filtration $\{\calF_t\}$ and satisfy the following properties:
	
	\noi
	(i) for every $0\le s<t$, $ B(t)- B(s)$ is independent of $\calF_s$ under $\PP$;
	
	\noi
	(ii) $ Y$ and $ R$ are nonnegative and nondecreasing;
	
	\noi
	(iii) 
	\begin{align}\label{313}
		X(t)= x_0+ mt+\sigma  B(t)+ Y(t)- R(t),\quad t\in\R_+,
	\end{align}
	such that 
	\begin{align}\notag
		X(t)\in [0,b],\quad t\in\R_+,\;\PP\text{-a.s.}
	\end{align}

	An {\it admissible control for the maximizer} is a meausre $ \Q$ defined on $(\Omega,\calF,\{\calF_t\})$ such that
	\begin{align}\label{315}
		\frac{d \Q(t)}{d\PP(t)}=\exp\Big\{\int_0^t \psi(s)dB(s)-\frac{1}{2}\int_0^t  \psi^2(s)ds\Big\},\quad t\in\R_+,
	\end{align}
	for an $\{\calF_t\}$-progressively measurable process $\psi$ satisfying 
	\begin{align}\label{316}
		&\E^{\PP}\Big[\int_0^\iy e^{-\varrho s} \psi^2(s)ds\Big]<\iy\quad\text{and}\quad\E^{\PP}\Big[e^{\frac{1}{2}\int_0^t \psi^2(s)ds}\Big]<\iy\quad\text{ for every $t\in\R_+$.}
	\end{align}
\end{definition}

The statements given in Remark \ref{rem_31} also hold for the RSDG as well. For completeness of the presentation and for later references we 
provide an alternative form of the dynamics given in \eqref{313},
\begin{align}\label{317}
	X(t)= x_0+ mt+\int_0^t\sigma\psi(s)ds+ \sigma  B^{\Q}(t)+ Y(t)- R(t),\quad t\in\R_+,
\end{align}
where $ B^{ \Q}(t):= B(t)-\int_0^t\psi(s)ds$, $t\in\R_+$, is an $\{\calF_t\}$-one-dimensional standard Brownian motion under $\Q$.

Denote by $ \calA( x_0)$ the set of all admissible controls for the minimizer, given the initial condition $ x_0$. As before, we often abuse notation and denote $ ( Y, R)\in \calA( x_0)$, keeping in mind that the control includes a filtered probability space. The set of all admissible controls for the maximizer is denoted by $ \calQ( x_0)$. 

{\bf The cost function (RSDG).} 
The expected cost associated with the initial condition $x$ and the controls $( Y, R)$ and $ \Q$ is given by
\begin{align}\notag
	J( x_0, Y, R, \Q;\eps):=&
	\E^{ \Q}\Big[\int_0^\iy e^{-\varrho t}\left( h(  X(t))dt + rd R(t)\right) \Big]-\frac{1}{\eps}L^\varrho( \Q\|\PP),
\end{align}
where 
\begin{align}\label{319a}
	&h(x):=\min\{\hat h\cdot\hat\xi : \hat\xi\in\calX,\; \theta\cdot\hat\xi=x\},\\\label{319b}
	&r:=\min\{\hat r\cdot q : q\in\R^I_+,\; \theta\cdot q=1\},
\end{align}
and $L^\varrho( \Q\|\PP)$ is given by \eqref{307a} with $(\Q,\PP)$ replacing $(\hat\Q_i,\PP_i)$. 
By the convexity of $\calX$ it follows that $h$ is convex. In fact, $h$ is piecewise linear and Lipschitz continuous. Moreover, $h(x)\ge 0$ for $x\ge 0$ and equality holds if and only if $x=0$. Therefore, $h$ is strictly increasing. In \cite[page 568]{ata-shi} it is shown that there is $i^*\in[I]$ such that,
\begin{align}\label{320}
	r=r_{i^*}\mu_{i^*}:=\min\{r_i\mu_i : i\in[I]\}.
\end{align}
The index $i^*$ stands for the class with the smallest rejection cost, weighted with the mean service rate. In fact, as we discuss in Section \ref{sec_33} and prove in Theorem \ref{thm_42}, under optimality of both players in the MSDG, rejections are performed only from this class. 

By the same arguments that lead to \eqref{308},
the cost function of the RSDG can alternatively be expressed by the technically more convenient form,
\begin{align}\label{321}
	J( x_0, Y, R, \Q;\eps)=&
	\E^{ \Q}\Big[\int_0^\iy e^{-\varrho t}\Big( h( X(t))dt + r  d R(t)-\frac{1}{2\eps} \psi^2(t)dt\Big) \Big],
\end{align}
with $ \psi$ satisfying \eqref{316} above. 
The value function is given by
\begin{align}\label{322}
	V( x_0;\eps)=\inf_{( Y, R)\in \calA( x_0)}\;\sup_{ \Q\in \calQ( x_0)}\;J( x_0, Y, R, \Q;\eps).
\end{align}

\begin{remark}\label{rem_32}
	In case that there is no ambiguity, we define the cost and the value functions by 
	\begin{align}\label{no_amb}
		J_{NA}( x_0, Y, R)&:=
		\E^{ \PP}\Big[\int_0^\iy e^{-\varrho t}( h( X(t))dt + rd R(t)) \Big],\\\notag
		V(x_0;0)&:=\inf_{(Y,R)\in\calA(x_0)}J_{NA}( x_0, Y, R).
	\end{align}
	This problem was studied by Harrison and Taksar \cite{har-tak} and later on was used by Atar and Shifrin \cite{ata-shi}. In Theorem \ref{thm_45} we show that this problem is obtained when the ambiguity vanishes, that is, $\lim_{\eps\to 0} V( x_0;\eps)=V(x_0;0)$.
\end{remark}

\subsection{The relationship between the games}\label{sec_33}
\begin{figure}
	\centering
	\includegraphics[width=0.5\textwidth]{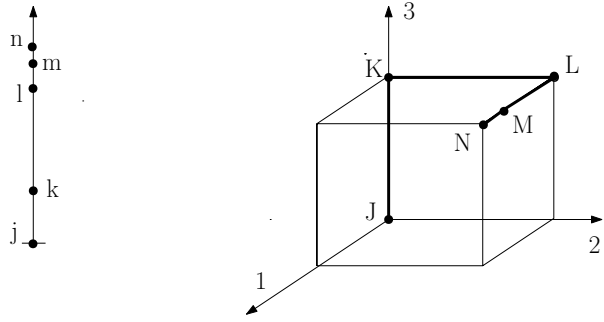}
	\caption{{\footnotesize The graphs refer to the case $I=3$, $\hat h=(1,5/2,3/2)$, $\mu=(3,1,3/2)$, and $(\hat b_1,\hat b_2,\hat b_3)=(4,7,6)$. The graph to the left stands for the workload levels. The curve of the function $\gamma$ is in bold in the graph to the right. The workload levels with the lower case letters are
			$j=0, k=\hat b_3/\mu_3=4, l=\hat b_3/\mu_3+\hat b_2/\mu_2=11, m= \hat b_3/\mu_3+\hat b_2/\mu_2+1$, and $n=\hat b_3/\mu_3+\hat b_2/\mu_2+\hat b_1/\mu_1=b=37/3$. They respectively correspond to the upper case letters: $J=(0,0,0), K=(0,0,\hat b_3/\mu_3)=(0,0,4), L=(0,\hat b_2/\mu_2,\hat b_3/\mu_3)=(0,7,4), M=(1,\hat b_2/\mu_2,\hat b_3/\mu_3)=(1,7,4)$, and $N=(\hat b_1/\mu_1,\hat b_2/\mu_2,\hat b_3/\mu_3)=(4/3,7,4)$.
		}\label{fig_1}}
\end{figure}

We now show that the last two games share the same value and moreover, that given any admissible control for the minimizer in either one of the games, one can construct an admissible control in the other game that performs at least as well. 
For this, we define a function $\gamma$, taken from \cite[Equations (48)--(49)]{ata-shi}, that sends any workload value $x_0$ to the cheapest state of the MSDG (from the holding cost perspective) among all the states whose workload levels are $x_0$. 
Using this function and an optimal strategy for the minimizer in the RSDG, we construct an optimal strategy for the minimizer in the MSDG, see Theorem \ref{thm_42}. In order to define the function it is convenient to assume without loss of generality that 
$$h_1\mu_1\ge h_2\mu_2\ge\dots\ge h_I\mu_I.$$
Recall that $b=\theta\cdot(\hat b_1,\ldots,\hat b_I)$. 
Given $x\in[0,b)$, let $(j,\upsilon)$ be the unique pair that is determined by
$$x=\sum_{i=j+1}^I\theta_i\hat b_i+\theta_j\upsilon,\qquad j\in[I],\quad\upsilon\in[0,\hat b_j),$$
and for $x=b$, take $(j,\upsilon)=(1,\hat b_1)$.
Let $\gamma:[0,b]\to\calX$ be the function given by
\begin{align}\label{323}
	\gamma(x)=\sum_{i=j+1}^I\hat b_ie_i+\upsilon e_j,
\end{align}
where $\{e_1,\ldots,e_I\}$ is the standard basis of $\R^I$. By filling up the cheaper buffers first (w.r.t.~the holding costs) we get that,
\begin{align}\label{newB}
	\gamma(x)\in\argmin\{\hat h\cdot\hat\xi:\hat\xi\in\calX,\theta\cdot\hat\xi=x\}.
\end{align}
The curve $\gamma(x)$, $x\in[0,b]$ is continuous and located on the edges of $\calX$, see Figure \ref{fig_1}. The idea is as follows, recall that the components of $\hat Y=(\hat Y_i)_{i=1}^I$ can be positive or negative, as long as $\theta\cdot\hat Y$ is nonnegative and nondecreasing. Now, as the workload changes in the interval $[0,b)$, the DM can use only the process $\hat Y$, without the need of the rejection process $\hat R$, so that $\hat X$ moves along the curve of $\gamma$. 
As will be shown in Theorem \ref{thm_42}, under optimality, the rejection process is used only to reduce the workload, and 
only from the class which has the cheapest weighted rejections cost, denoted by $i^*$. We discuss more about the minimizer's optimal strategy in the MSDG in Remark \ref{rem_33} and in the paragraph that comes before Theorem \ref{thm_42}.

We state the proposition below for an arbitrary initial point $\hat x_0\in\calX$ and show that $\hat V(\hat x_0;\kaboom)=V(x_0;\eps)$, where recall that $x_0=\theta\cdot\hat x_0$.

\begin{proposition}\label{prop_31} Fix $\kaboom=(\kaboom_i)_{i=1}^I$ and $\hat x_0\in\calX$. Let $\eps$ be given by \eqref{310g}.
	
	\noi
	(i) Given an admissible control in the MSDG, $(\Omega,\calF,\{\calF_t\},\PP,\hat B,\hat Y,\hat R)$ 
	and an admissible measure $\Q_*$ associated with $\psi_*$ in the RSDG that satisfies \eqref{315}--\eqref{316}, set $B=\frac{1}{\sigma}\theta\hat \sigma\cdot\hat B$,  $(X,Y,R)=(\theta\cdot\hat X,\theta\cdot\hat Y,\theta\cdot\hat R)$,
	and $\hat\psi_*=(\hat\psi_{*,1},\ldots,\hat\psi_{*,I})$, by
	\begin{align}\label{328}
		\hat\psi_{*,i}(t):=\frac{\sigma\psi_*(t)(\theta\hat\sigma)_i\kaboom_i}{\sum_{j=1}^I(\theta\hat\sigma)^2_j\kaboom_j}
		,\qquad t\in\R_+,\;i\in[I].
	\end{align}
	Let $\hat\Q_*$ be the associated measure defined through \eqref{304}. 
	Then $(Y,R) \in\calA(x_0)$, $\hat\Q_*\in\hat\calQ(\hat x_0)$, and 
	\begin{align}\label{324a}
		J(x_0,Y,R,\Q_*;\eps)\le\hat J(\hat x_0,\hat Y,\hat R,\hat\Q_*;\kaboom).
	\end{align}
	As a consequence, 
	\begin{align}\notag
		\sup_{\Q\in\calQ(x_0)}J(x_0,Y,R,\Q;\eps)\le\sup_{\hat\Q\in\hat\calQ(\hat x_0)}\hat J(\hat x_0,\hat Y,\hat R,\hat\Q;\kaboom).
	\end{align}
	\noi
	(ii) Conversely, consider an admissible control in the RSDG, $(\Omega,\calF,\{\calF_t\},\PP,B,Y,R)$, which is assumed to support an $I$-dimensional standard Brownian motion $\hat B$. 
	Consider also an admissible
	$\hat\Q_\sharp$ associated with $\hat\psi_\sharp$  in the MSDG that satisfies \eqref{304}--\eqref{305}.
	Assume that $\hat B$ is $\{\calF_t\}$-adapted and satisfies $\frac{1}{\sigma}\theta\hat \sigma\cdot\hat B=B$ and \eqref{300}. Define $(\hat X,\hat Y,\hat R)$ by
	\begin{align}\label{325}
		\hat X(t):=\gamma(X(t)),\quad \hat R(t):=R(t)\mu_{i^*}e_{i^*},
	\end{align}
	and
	\begin{align}\label{326}
		\hat Y(t):=\hat X(t)-\hat x_0-\hat mt-\hat\sigma \hat B(t) +\hat R(t).
	\end{align}
	Also, let 
	\begin{align}\label{335}
		\psi_\sharp(t):=\frac{1}{\sigma}\theta\hat\sigma\cdot\hat \psi_\sharp(t)
		,\quad t\in\R_+,
	\end{align}
	and let $\Q_\sharp$ be the associated measure defined through \eqref{315}. 
	Then $(\hat Y,\hat R)\in\hat\calA(\hat x_0)$, $\Q_\sharp\in\calQ(x_0)$, and 
	\begin{align}\notag
		\hat J(\hat x_0,\hat Y,\hat R,\hat\Q_\sharp;\kaboom)\le J(x_0,Y,R,\Q_\sharp;\eps).
	\end{align}
	As a consequence,
	\begin{align}\label{327c}
		\sup_{\hat\Q\in\hat\calQ(\hat x_0)}\hat J(\hat x_0,\hat Y,\hat R,\hat\Q;\kaboom)\le \sup_{\Q\in\calQ(x_0)}J(x_0,Y,R,\Q;\eps).
	\end{align}
\end{proposition}
Notice that in part (i) we define the dynamics $(X,Y,R)$ in the RSDG from the dynamics $(\hat X, \hat Y, \hat R)$ in the MSDG by projecting in the $\theta$ direction. In part (ii) we consider an auxiliary $I$-dimensional standard Brownian motion for technical reasons. The proof that the probability space $(\Omega,\calF,\{\calF_t\},\PP,B,Y,R)$ supports such a Brownian motion is merely technical and is fully given in \cite[Proposition 2.1.(iii)]{ata-shi}. Therefore, it is omitted. Also, in \eqref{325}--\eqref{326} we construct first $\hat X$ from the one-dimensional process $X$ using the well-defined function $\gamma$ that maps $[0,b]$ to $\calX$. Then, we define $\hat R$ as a process with degeneracy in its $[I]\setminus \{i^*\}$ coordinates. Only then we define $\hat Y$ in such a way that \eqref{302} holds.
\begin{corollary}\label{cor_31} Fix $\kaboom=(\kaboom_i)_{i=1}^I$ and $\hat x_0\in\calX$. Let $\eps$ be given by \eqref{310g}. 
	Then, 
	$\hat V(\hat x_0;\kaboom)= V(x_0;\eps)$.
\end{corollary}

\begin{remark}\label{rem_33}
	In Theorem \ref{thm_41} we show that the minimizer in the RSDG has a simple optimal strategy, which uses minimal idleness and minimal amount of rejections in order to keep the workload in some subinterval $[0,\beta]\subseteq[0,b]$. In Theorem \ref{thm_42} we use the function $\gamma$ given in \eqref{323} and the relations given in \eqref{325}--\eqref{326} to construct an optimal strategy for the minimizer in the MSDG. An optimal strategy for the maximizer emerges from \eqref{328} and \eqref{335}.
\end{remark}

\skp\noi
{\bf Proof of Proposition \ref{prop_31}:} 
Some of the arguments in the proof are given in \cite[Proposition 2.1]{ata-shi}. However, the ones that involve changes of the probability measure are new. For completeness of the proof we provide all the details.  

(i) The proof that $(Y,R) \in\calA(x_0)$ is straightforward and therefore omitted. 
Notice that $\|\hat\psi_*(\cdot)\|\le|\psi_*(\cdot)|$ and therefore,  
$\hat\Q_*\in\hat\calQ(\hat x_0)$ follows since $\Q_*\in\calQ(x_0)$. 

We now show that the distribution of $X$ under $\Q_*$ is the same as under $\hat\Q_* $. Replacing $\hat \Q$ by $\hat\Q_*$ in equation \eqref{306} and multiplying its both sides by $\theta$ yield
\begin{align}\label{329}
	X(t)&=x_0+mt+ \int_0^t \theta\hat\sigma\hat\psi_*(s)ds+\theta\hat\sigma\hat B^{\hat\Q_*}(t)+Y(t)-R(t)\\\notag
	&=x_0+mt+ \int_0^t \sigma \psi_*(s)ds+\sigma B^{\Q_*}(t)+Y(t)-R(t).
\end{align}
The equality between the integrals follows by the definitions of $\hat\psi_*$ and $\sigma$. 
The equality between the Brownian motion terms follows since 
\begin{align}\notag
	\theta\hat\sigma\hat B^{\hat\Q_*}(t)=\theta\hat\sigma\hat B(t)-\int_0^t\theta\hat\sigma\hat \psi_*(s)ds
	=\sigma B(t)-\sigma\int_0^t\psi_*(s)ds=\sigma B^{\Q_*},
\end{align}
which in turn follows by using the definitions of $\hat B^{\hat\Q_*}$ and $B^{\Q_*}$ (\eqref{306} and \eqref{317}) together with the definition of $B$ from the proposition and once again the definition of $\hat\psi_*$.
Recall that $\hat B^{\hat\Q_*}$ is an $I$-dimensional standard Brownian motion under $\hat\Q_*$ and that $\theta\hat\sigma$ is a deterministic vector with norm $\sigma$. Then from the above, we get that under $\hat\Q_*$, the process $B^{\Q_*}$ is a one-dimensional standard Brownian motion. Since $B^{\Q_*}$ is a one-dimensional standard Brownian motion also under $\Q_*$, we get from \eqref{329} that $X$ admits the same distribution under $\hat\Q_*$ and under $\Q_*$.

Next, by the definitions of $h$ and $r$ (see \eqref{319a} and \eqref{319b}) it follows that
\begin{align}\label{331}
	h(X(t))&\le\hat h\cdot\hat X(t),\qquad t\in\R_+,\\\label{331b}
	\int_0^\iy e^{-\varrho t} r d R(t)&\le \int_0^\iy e^{-\varrho t}\hat r\cdot d\hat R(t).
\end{align}
Moreover, using the fact that the distribution of $X$ under $\hat\Q_*$ is the same as under $\Q_*$ and the equality 
\begin{align}\notag
	\sum_{i=1}^I\frac{1}{2\kaboom_i}\hat\psi_{*,i}^2(t)=\frac{1}{2\eps}\psi_*^2(t),\qquad t\in\R_+,
\end{align}
we get,
\begin{align}\notag
	J( x_0, Y, R, \Q_*;\eps)&=
	\E^{ \Q_*}\Big[\int_0^\iy e^{-\varrho t}\Big( h( X(t))dt + r  d R(t)-\frac{1}{2\eps}\psi_*^2(t)\Big)dt \Big]\\\notag
	&=\E^{ \hat\Q_*}\Big[\int_0^\iy e^{-\varrho t}\Big( h( X(t))dt + r  d R(t)-\frac{1}{2\eps}\psi_*^2(t)\Big)dt \Big]\\\notag
	&\le \E^{ \hat\Q_*}\Big[\int_0^\iy e^{-\varrho t}\Big( \hat h\cdot \hat X(t)dt + \hat r \cdot d \hat R(t)-\sum_{i=1}^I\frac{1}{2\kaboom_i}\hat\psi_{*,i}^2(t)\Big)dt \Big]\\\notag
	&= \hat J(\hat x_0,\hat Y,\hat R,\hat\Q_*;\kaboom).
\end{align}
Since $\Q_*\in\calQ(x_0)$ is arbitrary, it follows that
\begin{align}\notag
	\sup_{\Q\in\calQ(x_0)}J(x_0,Y,R,\Q;\eps)\le\sup_{\hat\Q\in\hat\calQ(\hat x_0)}\hat J(\hat x_0,\hat Y,\hat R,\hat\Q;\eps).
\end{align}

(ii) We start by showing that $(\hat Y,\hat R)\in\hat\calA(\hat x_0)$. The processes $\hat X$, $\hat R$, and $\hat Y$ are $\{\calF_t\}$-adapted since $X$, $R$, and $\hat B$ are. 
Now, 
\begin{align}\notag
	\theta\cdot\hat Y(t) = \theta\cdot \hat X(t)-\theta\cdot x_0-\theta\cdot \hat mt-\theta\cdot \hat\sigma\hat B(t)+\theta\cdot \hat R(t)=Y(t),
\end{align}
and since $Y$ is admissible it follows that $\theta\cdot\hat Y=Y$ is nonnegative and nondecreasing. 
The processes $\{\hat R_i\}_{i=1}^I$ are clearly nonnegative and nondecreasing since so is $R$. Thus, property \eqref{301} holds. Finally, properties \eqref{302} and \eqref{303} follow by the definition of $\gamma$ and the construction of $(\hat X,\hat Y,\hat R)$. 

Notice that, $\|\hat\psi_\sharp(\cdot)\|=|\psi_\sharp(\cdot)|$ and therefore,  
$\Q_\sharp\in\calQ(x_0)$ follows since $\hat\Q_\sharp\in\hat\calQ(\hat x_0)$. 
From Cauchy-Schwartz inequality,
\begin{align}\notag
	\sum_{i=1}^I\frac{1}{2\kaboom_i}\hat\psi_{\sharp,i}^2(t)\ge\frac{1}{2\eps}\psi_\sharp^2(t).
\end{align}
By the definitions of $h$, $r$, and $\gamma$ and \eqref{320}, \eqref{newB}, and \eqref{325} it follows that
\begin{align}\notag
	\hat h\cdot\hat X(t) &=h(X(t)),\qquad t\in\R_+,\\\notag
	\int_0^\iy e^{-\varrho t}\hat r\cdot d\hat R(t)&= \int_0^\iy e^{-\varrho t} r d (\theta\cdot\hat R(t))= \int_0^\iy e^{-\varrho t} r d R(t).
\end{align}
By the relationship between the $I$-dimensional processes $\hat B$ and $\hat\psi_\sharp$ and the one-dimensional processes $B$ and $\psi_\sharp$, one can show as was done in the previous part that the distribution of $X$ under $\hat\Q_\sharp$ is the same as under $\Q_\sharp$. Therefore,
\begin{align}\label{337}
	J( x_0, Y, R, \Q_\sharp;\eps)&=
	\E^{ \Q_\sharp}\Big[\int_0^\iy e^{-\varrho t}\Big( h( X(t))dt + r  d R(t)-\frac{1}{2\eps} \psi_\sharp^2(t)\Big)dt \Big]\\\notag
	&=\E^{ \hat\Q_\sharp}\Big[\int_0^\iy e^{-\varrho t}\Big( h( X(t))dt + r  d R(t)-\frac{1}{2\eps} \psi_\sharp^2(t)\Big)dt \Big]\\\notag
	&\ge \E^{ \hat\Q_\sharp}\Big[\int_0^\iy e^{-\varrho t}\Big( \hat h\cdot \hat X(t)dt + \hat r \cdot d \hat R(t)-\sum_{i=1}^I\frac{1}{2\kaboom_i}\hat\psi_{\sharp,i}^2(t)\Big)dt \Big]\\\notag
	&= \hat J(\hat x_0,\hat Y,\hat R,\hat\Q_\sharp;\kaboom).
\end{align}
Since $\hat\Q_\sharp\in\hat\calQ(\hat x_0)$ is arbitrary, it follows that
\begin{align}\notag
	\sup_{\hat\Q\in\hat\calQ(\hat x_0)}\hat J(\hat x_0,\hat Y,\hat R,\hat\Q;\kaboom)\le\sup_{\Q\in\calQ(x_0)}J(x_0,Y,R,\Q;\eps).
\end{align}

\hfill$\Box$

\section{Solution of the RSDG}\label{sec_4}
\beginsec
In this section we provide a solution to the RSDG. In Section \ref{sec_41a} we present the notion of a reflection strategy. Then in Section \ref{sec_41} we provide the HJB equation associated with the RSDG. We prove that the value function is the unique smooth solution of the HJB equation and show that the minimizer has an optimal reflecting strategy. 

\subsection{Reflecting strategies}\label{sec_41a}
The optimal strategy of the minimizer is shown to be one that enforces the workload to stay in a specific interval of the form $[0,\beta]$ with minimal effort. To rigorously define such a strategy we make use of the {\it Skorokhod map on an interval}. Fix $\al<\beta$. For any $\eta\in\calD(\R_+,\R)$ there exists a unique triplet of functions $(\chi,\zeta_1,\zeta_2)\in\calD(\R_+,\R^3)$ that satisfies the following properties:

\noi
(i) for every $t\in\R_+$, $\chi(t)=\eta(t)+\zeta_1(t)-\zeta_2(t)$;

\noi
(ii) $\zeta_1$ and $\zeta_2$ are nondecreasing, $\zeta_1(0-)=\zeta_2(0-)=0$, and
\begin{align}\notag
	\int_0^\iy
	1_{(\al,\beta]}(\chi(t))d\zeta_1(t)=\int_0^\iy
	1_{[\al,\beta)}(\chi(t))d\zeta_2(t)=0.
\end{align}
We denote by $\Gamma_{[\al,\beta]}(\eta)=(\Gamma_{[\al,\beta]}^1,\Gamma_{[\al,\beta]}^2,\Gamma_{[\al,\beta]}^3)(\eta)=(\chi,\zeta_1,\zeta_2)$. 
See \cite{Kruk2007} for existence and uniqueness of solutions,
and continuity and further properties of the map.
In particular, we have the following.
\begin{lemma}\label{lem_Skorokhod}
	There exists a constant $c_S>0$ such that
	for every $t>0$, $\al<\beta$ and $\om,\tilde\om\in\calD(\R_+,\R)$,
	\begin{equation}\notag
		\sup_{s\in[0,t]}\|\Gam_{[\al,\beta]}(\om)(t) - \Gam_{[\al,\beta]}(\tilde\om)(t)\|
		\le c_S\sup_{s\in[0,t]}|\om(s)-\tilde\om(s)|.
	\end{equation}
\end{lemma}

\begin{definition}\label{def_Skorokhod}
	Fix $x_0,\beta\in[0,b]$. 
	The strategy $(Y,R)$ is called 
	a $\beta$-reflecting strategy if for every $\eta\in \calC(\R_+,\R)$ one has $(X,Y,R)(\eta)=\Gamma_{[0,\beta]}(\eta)$. 
\end{definition}
One can easily verify that any $\beta$-reflecting strategy is admissible.

\subsection{The HJB equation and the value function}\label{sec_41}
In case that there is no ambiguity (see Remark \ref{rem_32}), the problem was analyzed by Harrison and Taksar \cite{har-tak}. 
It is shown there (see Proposition 5.11) that if $k_{\al\beta}$ is a $\calC^1([0,b],\R)$ is twice continuously differentiable on $[\al,\beta]$ and satisfies
\begin{align}\label{402a}
	l+k_{\al\beta}'(x)&=0,\qquad  0\le x\le\al,\\\label{402b}
	\frac{1}{2}\sigma^2 k_{\al\beta}''(x)+mk_{\al\beta}'(x)-\varrho k_{\al\beta}(x)+h(x)&=0,\qquad \al\le x\le\beta,\\\label{402c}
	r-k_{\al\beta}'(x)&=0,\qquad \beta\le x\le b,
\end{align}
then 
\begin{align}\notag
	k_{\al\beta}(x)=\E^\PP\left[\int_0^\iy e^{-\varrho t}[h(X(t))dt+rdR(t)+ldY(t)]\right], \quad x\in[0,b],
\end{align}
with $(X,Y,R)(t)=\Gamma_{[\al,\beta]}(x_0+m\cdot+\sigma B(\cdot))(t)$, $t\in\R_+$.
The rationale behind this argument is as follows. Consider $x\in(\beta,b)$, then in order to keep the process $X$ between $\al$ and $\beta$, there is an instantaneous reflection from above that contributes the cost $r(x-\beta)$. This explains \eqref{402c}. Similar arguments yield \eqref{402a}. When no control is taking action, standard arguments imply \eqref{402b}.
The HJB in this case takes the form
\begin{align}\notag
	\begin{cases} \Big[\frac{1}{2}\sigma^2 f''(x)+mf'(x)-\varrho f(x)+h(x)\Big]\wedge f'(x) \wedge [r-f'(x)]=0, &x\in (0,b) ,\\ 
		f'(0)=0,\quad f'(b)=r.
	\end{cases} 
\end{align}
However, the uniqueness of a solution of the HJB is not argued in \cite{har-tak} (see the paragraph before Section 7) but rather in \cite[Proposition 2.2]{ata-shi}, using viscosity solutions. 

\begin{remark}\label{rem_41}
	Notice that in \cite{har-tak} there is a cost associated with the reflections from both sides, unlike our case that does not consider a cost component for reflections from below, that is $l=0$. It simply follows since in our QCP there is no penalty/reward for idleness. Henceforth, under optimality, $\al=0$ in \eqref{402a}--\eqref{402b}. 
\end{remark}
Motivated by these results and the structure of the cost function given in \eqref{321} together with the inf-sup structure of the value function given in \eqref{322}, we consider the following HJB, 
\begin{align}\notag
	\begin{cases} 
		\Big[\underset{p\in\R}{\sup}\{\frac{1}{2}\sigma^2 f''(x)+(m+\sigma p)f'(x)-\varrho f(x)+h(x)-\frac{1}{2\eps}p^2\}\Big]\wedge f'(x) \wedge [r-f'(x)]=0, &x\in(0,b), \\ 
		f'(0)=0,\quad f'(b)=r.
	\end{cases} 
\end{align}
Or equivalently, by substituting the optimal solution of the $\sup_{p\in\R}$ above, $p^*=\eps\sigma f'(x)$,
\begin{align}\tag{HJB($\eps$) }
	\begin{cases} 
		[f''(x)+H(x,f(x),f'(x))]\wedge f'(x) \wedge [r-f'(x)]=0, &x\in(0,b), \\ 
		f'(0)=0,\quad f'(b)=r,
	\end{cases} 
\end{align}
where hereafter,
\begin{align}\notag
	H(x,y,z):=\frac{2}{\sigma^2} \left(mz+\frac{1}{2}\sigma^2\eps z^2-\varrho y+h(x)\right).
\end{align}
Notice that when $\eps=0$, HJB($\eps$) coincides with the HJB given in Harrison and Taksar \cite[Equation (1.2)]{har-tak} and by Atar and Shifrin \cite[Equation (41)]{ata-shi}.

The idea is that we may think of the change of measure done by the maximizer as a change of the drift term from $m$ to $m+\sigma p$, which costs her  $\tfrac{1}{2\eps}p^2$, where $p$ depends on $x$ and is chosen in order to maximize the cost. As can be seen, due to the term $\frac{1}{2}\sigma^2\eps(f'(x))^2$ in $H(x,f(x),f'(x))$, the HJB is not linear, a fact that raises some technical difficulties in proving Proposition \ref{prop_43} below.

We now state the main theorem of the paper, which is given also for $\eps=0$, see Remark \ref{rem_32}.

\begin{theorem}\label{thm_41}
	Fix $\eps\in[0,\iy)$. The value function $V(\cdot;\eps)$ is the unique $\calC^2([0,b],\R)$ solution of HJB($\eps$). Moreover, set 
	\begin{align}\label{440}
		\beta_\eps=\inf\left\{x\in(0,b] : V'(x;\eps)=r\right\},
	\end{align}
	where $V'(x;\eps)$ is the derivative of $V(\cdot;\eps)$ w.r.t.~$x$. 
	Then the $\beta_\eps$-reflecting strategy is optimal for the minimizer and $V=V(\cdot;\eps)$ satisfies,
	\begin{align}\label{414z}
		\begin{cases}
			V''(x)+H(x,V(x),V'(x))=0,\qquad &0\le x\le\beta_\eps,\\
			r-V'(x)=0,\qquad &\beta_\eps\le x\le b,\\
			V'(0)=0.
		\end{cases}
	\end{align}
\end{theorem}
From the definition of HJB($\eps$) and the theorem above we get the following corollary, which is given for reference purposes.
\begin{corollary}\label{cor_42b}
	For any $\eps\in(0,\iy)$, $0\le V'(\cdot;\eps)\le r$.
\end{corollary}

%
%
%
The proof of the theorem is done in several steps and follows from the next three propositions. 
Before stating them we define a measure $\Q_f$ associated with a process $\psi_f$, which in turn is driven by a function $f\in\calC^1([0,b],\R)$ through \eqref{315}. This measure serves us in the sequel, especially in Propositions \ref{prop_41} and \ref{prop_42} and Theorem \ref{thm_42}. In the latter, we show that the measure $\Q_V=\Q_{V(\cdot;\eps)}$ is the optimal strategy of the maximizer, where $V$ is the value function.
For any $f\in\calC^1([0,b],\R)$ and $t\in\R_+$ set
\begin{align}\label{407}
	\psi_f(t):&=\argmax_{p\in\R}\left\{\frac{1}{2}\sigma^2f''(X(t))+(m+\sigma p)f'(X(t))-\varrho f(X(t))+h(X(t))-\frac{1}{2\eps}p^2\right\}\\\notag
	&=
	\eps\sigma f'(X(t)).  
\end{align}
This is an $\{\calF_t\}$-progressively measurable process since $Y,R$, and $B$ are, see Definition \ref{def_33}. Let $\Q_f$ be the measure associated with $\psi_f$ through \eqref{315}. Note that $\psi_f$ is bounded since $f\in\calC^1([0,b],\R)$ and therefore \eqref{316} holds trivially. 

The proof of the propositions below are deferred to after the proof of Theorem \ref{thm_41}. 
\begin{proposition}\label{prop_41}
	Fix $\eps\in(0,\iy)$. Assume that HJB($\eps$) admits a $\calC^2([0,b],\R)$ solution $f$. Then 
	\begin{align}\notag
		f(x)\le \inf_{(Y,R)\in\calA(x)} J(x,Y,R,\Q_f;\eps),\qquad x\in[0,b],
	\end{align}
	and as a consequence $f\le V$. 
\end{proposition}  
The next proposition characterizes the solution of the following ordinary differential equation, 
\begin{align}\label{414}
	\begin{cases}
		k_\beta''(x)+H(x,k_\beta(x),k_\beta'(x))=0,\qquad &0\le x\le\beta,\\
		r-k_{\beta}'(x)=0,\qquad &\beta\le x\le b,\\
		k_{\beta}'(0)=0.
	\end{cases}
\end{align}

\begin{proposition}\label{prop_42}
	Fix $\eps\in(0,\iy)$. Assume that there is a function $k_{\beta}=k_{\beta,\eps}\in\calC^1([0,b],\R)\cap\calC^2([a,b]\setminus\{\beta\},\R)$ that solves \eqref{414}. Let $(Y_\beta,R_\beta)$ be a $\beta$-reflecting strategy, that is $(X,Y_\beta,R_\beta)(t)=\Gamma_{[0,\beta]}(x+m\cdot+\sigma B(\cdot))(t)$, $t\in\R_+$. Then, for every $x\in[0,b]$,
	\begin{align}\label{415}
		k_{\beta}(x)&=\sup_{\Q\in\calQ(x)}J(x,Y_\beta,R_\beta,\Q;\eps)=J(x,Y_\beta,R_\beta,\Q_k;\eps),
	\end{align}
	where $\Q_k=\Q_{k_\beta}$  is the measure associated with $\psi_k=\psi_{k_\beta}$, defined in \eqref{407}. 
\end{proposition}
We now claim that for every $\eps\in(0,\iy)$, HJB($\eps$) admits a unique smooth solution. Harrison and Taksar provided in \cite[page 450]{har-tak} explicit functions from which a smooth solution can be constructed. The construction is provided and the smoothness is claimed to be straightforward yet tedious and therefore omitted, see the paragraph below (6.8) there. The situation is more subtle in our case since HJB($\eps$) is nonlinear and therefore explicit solutions are out of reach. We choose a different path and use the {\it shooting method} to prove that a smooth solution of HJB($\eps$) uniquely exists and that it solves the free-boundary problem \eqref{414}. In short, the shooting method is used to solve boundary value problems by reducing them to initial value problems; see \cite[Section 7.3]{Stoer1980} for further reading about the method. We take it one step forward and use it in the free-boundary setup with Neumann boundary conditions. Specifically, we consider a set of Cauchy problems on the interval $[0,b]$, indexed\footnote{For technical reasons, in the proof we use a modification of $H$.} by $s$: 
\begin{align}\notag
	\begin{cases}
		(k^{(s)})''(x)+H(x,k^{(s)}(x),(k^{(s)})'(x))=0,\qquad &0\le x\le b,\\
		k^{(s)}(0)=0,\quad(k^{(s)})'(0)=0.
	\end{cases}
\end{align}
We prove that there is a parameter $s$ for which there is $\beta_\eps\in(0,b]$ at which we can smoothly paste the linear functions $k^{(s)}$ and $x\mapsto k^{(s)}(\beta_\eps)+r(x-\beta_\eps)$; namely, $(k^{(s)})'(\beta_\eps)=r$ and in case that $\beta_\eps<b$ also $(k^{(s)})''(\beta_\eps)=0$. This implies a smooth solution to \eqref{414z}. 
From the technical point of view, the proof of Proposition \ref{prop_43} is the most demanding in this section.
\begin{proposition}\label{prop_43}
	For every $\eps\in(0,\iy)$, HJB($\eps$) admits a unique $\calC^2([0,b],\R)$ solution. Moreover, the solution has the form $k_{\beta_\eps}$ for some parameter $\beta_\eps\in(0,b]$ for which $(k_{\beta_\eps})'<r$ on $[0,\beta_\eps)$, where $k_{\beta_\eps}$ is taken from Proposition \ref{prop_42}.
\end{proposition}

\noi
{\bf Proof of Theorem \ref{thm_41}:} Recall that when $\eps=0$, HJB($\eps$) coincides with the HJB given in \cite{har-tak, ata-shi}. Thus, we focus only on positive $\eps$'s.

From Proposition \ref{prop_43}, the HJB($\eps$) admits a unique $\calC^2([0,b],\R)$ solution that also solves \eqref{414} for some $\beta_\eps\in(0,b]$. Denote it by $k_{\beta_\eps}$. From Proposition \ref{prop_42} this is the cost of the $\beta_\eps$-reflecting strategy. Therefore, the value function, which is the infimum over all the strategies satisfies, $V(\cdot;\eps)\le k_{\beta_\eps}(\cdot)$. Together with Proposition \ref{prop_41}, $V(\cdot;\eps)= k_{\beta_\eps}(\cdot)$. Recalling again that $k_{\beta_\eps}$ is the cost of the $\beta_\eps$-reflecting strategy, we obtain its optimality. The relation in \eqref{440} follows now from Proposition \ref{prop_43}.

\hfill$\Box$

\skp
The rest of the section is devoted to the proofs of Propositions \ref{prop_41}--\ref{prop_43}.

\noi
{\bf Proof of Proposition \ref{prop_41}:} Fix $f\in\calC^2([0,b],\R)$ and an arbitrary $(Y,R)\in\calA(x)$ with an associated standard Brownian motion $B$. Set
\begin{align}\notag
	X(t)= x+ mt+\sigma  B(t)+ Y(t)- R(t),\quad t\in\R_+.
\end{align}
%
%
%
Recalling \eqref{317}, It{\^o}'s lemma implies that for every $t>0$ and every $\Q\in\calQ(x)$,
\begin{align}\notag
	\E^{\Q}\left[e^{-\varrho t}f(X(t))\right]=f(x)&+\E^{\Q}\left[\int_0^te^{-\varrho s}\left(\frac{1}{2}\sigma^2f''(X(s))+(m+\sigma\psi(s))f'(X(s))-\varrho f(X(s))\right)ds\right]\\\label{408}
	&+\E^{\Q}\left[\int_0^te^{-\varrho s}f'(X(s))(dY^c(s)-dR^c(s))\right]\\\notag
	&+\E^{\Q}\left[\sum_{0\le s\le t}e^{-\varrho s}\Delta f(X)(s)\right],
\end{align}
where $\Q$ and $\psi$ are related to each other through \eqref{315}--\eqref{316}. We used the fact that $B^{\Q}(t):=B(t)-\int_0^t\psi(s)ds$, $t\in\R_+$ is a $\Q$ standard Brownian motion. The processes
\begin{align}\label{409}
	Y^c(t):=Y(t)-\sum_{0\le s\le t}\Delta Y(s),\qquad R^c(t):=R(t)-\sum_{0\le s\le t}\Delta R(s),\qquad  t\in\R_+,
\end{align}
are the continuous parts of $Y$ and $R$, respectively. Consider now \eqref{408} with $\psi=\psi_f$ defined in \eqref{407} and with the measure $\Q=\Q_f$. 
Recalling the definition of $\psi_f$ and that $f$ solves HJB($\eps$), we get that 
\begin{align}\label{410}
	\E^{\Q_f}\left[e^{-\varrho t}f(X(t))\right]\ge &f(x)-\E^{\Q_f}\left[\int_0^te^{-\varrho s}\left(h(X(s))ds+rdR(s)-\frac{1}{2\eps}\psi_f^2(s)ds\right)\right]\\\notag
	&+\E^{\Q_f}\left[\sum_{0\le s\le t}e^{-\varrho s}(\Delta f(X)(s)+r\Delta R(s))\right].
\end{align}
Since $\Delta X(s)=\Delta Y(s)-\Delta R(s)$, we get that
\begin{align}\label{411}
	\Delta f(X)(s)+r\Delta R(s)&=f(X(s))-f(X(s)-\Delta Y(s))\\\notag
	&\quad-[f(X(s)+\Delta R(s)-\Delta Y(s))-f(X(s)-\Delta Y(s))]+r\Delta R(s)\\\notag
	&=\int_{X(s)-\Delta Y(s)}^{X(s)}f'(u)du +\int_{X(s)-\Delta Y(s)}^{X(s)+\Delta R(s)-\Delta Y(s)}(r- f'(u))du \\\notag
	&\ge 0.
\end{align}
Combining \eqref{410}--\eqref{411}, we get that 
\begin{align}\notag
	\E^{\Q_f}\left[e^{-\varrho t}f(X(t))\right]+\E^{\Q_f}\left[\int_0^te^{-\varrho s}\left(h(X(s))ds+rdR(s)-\frac{1}{2\eps}\psi_f^2(s)ds\right)\right]\ge f(x).
\end{align}
Recalling the definition of the cost function $J$ in \eqref{321} and noting that the function $f$ is bounded as a continuous function on $[0,b]$, by taking $t\to\iy$, we get that 
\begin{align}\notag
	\sup_{\Q\in\calQ(x)}J(x,Y,R,\Q;\eps)\ge J(x,Y,R,\Q_f;\eps)\ge f(x).
\end{align}
Since $x$ and $(Y,R)$ are arbitrary, it follows that $V(\cdot;\eps)\ge f(\cdot)$.

\hfill$\Box$

\skp\noi{\bf Proof of Proposition \ref{prop_42}:}
%
We split the proof into two cases $x\in[0,\beta]$ and $x\in(\beta,b)$. Fix $x\in[0,\beta]$ and an arbitrary $\Q\in\calQ(x)$ with an associated process $\psi$. Recall the notation $\psi_k=\psi_{k_\beta}$. By \eqref{407}, we get
\begin{align}\label{415c}
	&\frac{1}{2}\sigma^2k_\beta''(X(t))+(m+\sigma \psi(t))k_\beta'(X(t))-\varrho k_\beta(X(t))+h(X(t))-\frac{1}{2\eps}(\psi(t))^2\\\notag
	&\quad\le\frac{1}{2}\sigma^2k_\beta''(X(t))+(m+\sigma \psi_k(t))k_\beta'(X(t))-\varrho k_\beta(X(t))+h(X(t))-\frac{1}{2\eps}(\psi_k(t))^2\\\notag
	&\quad=\frac{1}{2}\sigma^2\left(k_\beta''(x)+H(x,k_\beta(x),k_\beta'(x))\right)=0.
\end{align}

Equations \eqref{408} and \eqref{409} are given for general admissible strategies $(Y,R)\in\calA(x)$ and $\Q\in\calQ(x)$. Thus, they hold here as well. Notice that since $x\in[0,\beta]$ and $(Y_\beta,R_\beta)$ is a $\beta$-reflecting strategy, the processes $Y_\beta$ and $R_\beta$ have no jumps and so, $Y_\beta^c=Y_\beta$ and $R_\beta^c=R_\beta$. From \eqref{408}, \eqref{409}, and \eqref{415c}, one has,
\begin{align}\label{415d}
	\E^{\Q}\left[e^{-\varrho t}k_\beta(X(t))\right]\le k_\beta(x)&-\E^{\Q}\left[\int_0^te^{-\varrho s}\left(h(X(s)ds-\frac{1}{2\eps}\psi^2(s)ds\right)ds\right]\\\notag
	&+\E^{\Q}\left[\int_0^te^{-\varrho s}k_\beta'(X(s))(dY_\beta(s)-dR_\beta(s))\right].
\end{align}
Using now the equalities $k_\beta'(0)=0$ and $k_\beta'(\beta)=r$, driven from \eqref{414}, we get
\begin{align}\label{415e}
	\E^{\Q}\left[e^{-\varrho t}k_\beta(X(t))\right]\le k_\beta(x)&-\E^{\Q}\left[\int_0^te^{-\varrho s}\left(h(X(s))ds+rdR_\beta(s)-\frac{1}{2\eps}\psi^2(s)ds\right)\right].
\end{align}
Recalling the definition of $J$ and noting that the function $k_\beta$ is bounded as a continuous function on $[0,b]$, by taking $t\to\iy$, we get that 
\begin{align}\label{415f}
	k_\beta(x)\ge \sup_{\Q\in\calQ(x)}J(x,Y_\beta,R_\beta,\Q;\eps).
\end{align}
Notice that all the inequalities in \eqref{415c}--\eqref{415e} hold with equality in case that $\psi=\psi_k$ and $\Q=\Q_k$. Then, together with \eqref{415f}, we get that \eqref{415} holds.

Consider now the case that $\beta<b$ and $x\in(\beta,b]$. From \eqref{414}, $$k_{\beta}(x)=r(x-\beta)+k_{\beta}(\beta).$$ Since the strategy $(Y_\beta,R_\beta)$ starts with an instantaneous rejection of $x-\beta_0$, there is an immediate cost of $r(x-\beta)$ and hence for any $\Q\in\calQ(x)$,
$$J(x,Y_\beta,R_\beta,\Q;\eps)=r(x-\beta)+J(\beta,Y_\beta,R_\beta,\Q;\eps).$$ Recall that we just showed that 
$$k_{\beta}(\beta)=\sup_{\Q\in\calQ(x)}J(\beta,Y_\beta,R_\beta,\Q;\eps)=J(\beta,Y_\beta,R_\beta,\Q_f;\eps).$$ 
From the last three equalities we get that \eqref{415} holds for $x\in(\beta,b]$ as well.

\hfill$\Box$

\skp
Before getting to the proof of Proposition \ref{prop_42} we provide an auxiliary parametrized ordinary differential equation in order to solve the free-boundary one. Fix a parameter $s\in\R$ and consider the following Cauchy problem 
\begin{align}\label{423}
	\begin{cases}
		(k^{(s)})''(x)+H_F(x,k^{(s)}(x),(k^{(s)})'(x))=0,\qquad x\in[0,b],\\
		(k^{(s)})'(0)=0,\quad k^{(s)}(0)=s,
	\end{cases}
\end{align}
where 
\begin{align}\label{HF}
	H_F(x,y,z):=H(x,y,F(z))
\end{align}
and $F$ is a $\calC^1(\R,\R)$ function that satisfies the following properties: $F(z)=z$ on $[-r,r]$, $|F|\le 2r$, and $|F'|\le 1$ and therefore Lipschitz continuous. For example,
\begin{align}\notag
	F(z)=\begin{cases}
		-3/2r, &\qquad\quad z< -2r,\\
		(1/2)r+2z+z^2/(2r), & -2r\le z<-r,\\
		z, & \;\;-r\le z\le r, \\
		-(1/2)r+2z-z^2/(2r), & r\le z<2r,\\
		3/2r, &\;\;\; 2r\le z.
	\end{cases}
\end{align}
Because the function $F$ and its derivative are bounded, and since the function $h$ is Lipschitz (see the paragraph below \eqref{319b}), $H_F$ is uniformly Lipschitz. Namely, there is a constant $c_L$ such that for every $(x,y,z),(x',y',z')\in[0,b]\times\R\times\R$, one has
\begin{align}\label{ac05}
	|H_F(x,y,z)-H_F(x',y',z')|\le L(|x-x'|+|y-y'|+|z-z'|).
\end{align}
From \cite[Section 0.3.1]{Polyanin2003}, 
\eqref{423} admits a unique $\calC^2([0,b],\R)$ solution.   

Set
\begin{align}\label{ac09}
	\beta^{(s)}&:=\inf\{x\in(0,b] : (k^{(s)})'(x)\ge r\}\wedge b,
\end{align}
where we use the convention that $\inf\emptyset =\iy$. The smoothness of $k^{(s)}$ implies that 
\begin{align}\label{ac04}
	\text{if $\;\beta^{(s)}<b\;$ then $\;(k^{(s)})'(\beta^{(s)})=r$.}
\end{align}

The following lemma provides some continuity properties that serve us in the proof of Proposition \ref{prop_43}.
\begin{lemma}\label{lem_42}
	The function $s\mapsto(k^{(s)},(k^{(s)})',(k^{(s)})'')$ is continuous in the uniform norm topology taken on the interval $[0,b]$. Moreover, the mapping $s\mapsto\beta^{(s)}$ is continuous for every $s$ for which either $\beta^{(s)}=b$ or the following two conditions hold $\beta^{(s)}<b$ and $(k^{(s)})''(\beta^{(s)})\ne 0$. In these cases we conclude that the mapping  $s\mapsto(k^{(s)}(\beta^{(s)}),(k^{(s)})'(\beta^{(s)}),(k^{(s)})''(\beta^{(s)}))$ is also continuous. 
\end{lemma}
{\bf Proof:} 
Fix $s\in\R$. One can easily verify that the conditions of \cite[Theorem 23]{Protter1967} are satisfied for $H_F$. Therefore, 
for any $\delta_1\in\R$ 
\begin{align}\label{432}
	\sup_{x\in[0,b]}\left|k^{(s+\delta_1)}(x)-k^{(s)}(x)\right|\le|\delta_1|,
\end{align}
and the continuity of $s\mapsto k^{(s)}$ is established. 

We now turn to showing that $s\mapsto (k^{(s)})'$ is continuous. From \eqref{423} it follows that for every $x\in[0,b]$
\begin{align}\notag
	(k^{(s)})'(x)&=0-\int_0^x H_F(y, k^{(s)}(y), (k^{(s)})'(y))dy,\\\notag
	(k^{(s+\delta_1)})'(x)&=0-\int_0^x H_F(y, k^{(s+\delta_1)}(y), (k^{(s+\delta_1)})'(y))dy.
\end{align}
Set $g_{s,\delta_1}(x):=\left| (k^{(s+\delta_1)})'(x)- (k^{(s)})'(x)\right|$, $x\in[0,b]$. From \eqref{ac05}, there exists a constant $L>0$ independent of $s$ and $\delta_1$ such that 
\begin{align}\notag
	g_{s,\delta_1}(x)&\le L\int_0^x \left(\left|k^{(s+\delta_1)}(y)-k^{(s)}(y)\right|+g_{s,\delta_1}(y)\right)dy
	\le Lb|\delta_1|+L\int_0^x g_{s,\delta_1}(y)dy,
\end{align}
where the second inequality follows by \eqref{432}. Now, Gr\"{o}nwall's inequality implies that \\$\sup_{x\in[0,b]}g_{s,\delta_1}(x)\le |\delta_1| Lbe^{Lb}$, and the continuity of $ s\mapsto(k^{(s)})'$ is established.

Finally, the continuity of $ s\mapsto(k^{(s)})''$ follows by the relation $(k^{(s)})''(x)=-H_F(x,k^{(s)}(x),(k^{(s)})'(x))$, the continuity of $ s\mapsto(k^{(s)},(k^{(s)})')$, and the Lipschitz continuity of $H_F$ stated in \eqref{ac05}.

The rest of the proof is dedicated to the continuity of the function $s\mapsto \beta^{(s)}$ under the conditions mentioned in the lemma. To this end, we fix $s\in\R$ and show that if $\beta^{(s)}=b$ or if $\beta^{(s)}<b$ and $(k^{(s)})'(\beta^{(s)})< r$, then
\begin{align}\label{435}
	\limsup_{\delta\to 0}\;\beta^{(s+\delta)}\le\beta^{(s)}\le\liminf_{\delta\to 0}\;\beta^{(s+\delta)}.
\end{align}
We start with the first inequality. If $\beta^{(s)}=b$ then it is obvious, since all the $\beta^{(u)}$'s are less or equal to $b$. If $\beta^{(s)}<b$ and $(k^{(s)})''(\beta^{(s)})\ne 0$, we necessarily have $(k^{(s)})''(\beta^{(s)})> 0$. Otherwise, $(k^{(s)})''(\beta^{(s)})< 0$ and from \eqref{ac04}, we get that $(k^{(s)})'(\beta^{(s)}-\nu)> r$ for sufficiently small $\nu>0$, a contradiction to the definition of $\beta^{(s)}$.

Using now $(k^{(s)})''(\beta^{(s)})> 0$, we get that for sufficiently small $\nu>0$, 
$(k^{(s)})'(\beta^{(s)}+\nu)>r$. By the continuity of $s\mapsto (k^{(s)})'$, we get that for every $\delta_2$ with sufficiently small absolute value, one has $(k^{(s+\delta_2)})'(\beta^{(s)}+\nu)>r$. Therefore, $\beta^{(s+\delta_2)}<\beta^{(s)}+\nu$ and
$\limsup_{\delta\to 0}\;\beta^{(s+\delta)}\le\beta^{(s)}+\nu$. Since $\nu>0$ can be arbitrary small we get the first inequality on \eqref{435}.

We now turn to proving the second inequality in \eqref{435}. Set $\gamma_1>0$ and $\hat\beta^{(s)}:=\liminf_{\delta\to 0}\;\beta^{(s+\delta)}$. Consider a sequence $\{\delta_j\}_j\to0$ such that $\{\beta^{(s+\delta_j)}\}_j\to \hat\beta^{(s)}$ and $\beta^{(s+\delta_j)}<b$ for every $j$. If such a subsequence does not exist it means that for every $\delta$ with sufficiently small absolute value, $\beta^{(s+\delta)}=b$ and the claim is trivial. Since $k^{(s)}\in\calC^2([0,b],\R)$, we get that for sufficiently large $j$,
\begin{align}\label{436}
	\left|(k^{(s)})'(\beta^{(s+\delta_j)})-(k^{(s)})'(\hat \beta^{(s)})\right|<\gamma_1.
\end{align}
Together with the continuity of 
$s\mapsto (k^{(s)})'$, we get that for sufficiently large $j$, one has
\begin{align}\label{437}
	\left|(k^{(s+\delta_j)})'(\beta^{(s+\delta_j)})-(k^{(s)})'( \beta^{(s+\delta_j)})\right|<\gamma_1.
\end{align}
Recall that $\beta^{(s+\delta_j)}<b$. Thus, $(k^{(s+\delta_j)})'(\beta^{(s+\delta_j)})=r$.
From \eqref{436}--\eqref{437} we get that 
\begin{align}\notag
	\left|(k^{(s)})'(\hat \beta^{(s)})-r\right|<2\gamma_1.
\end{align}
Since $\gamma_1>0$ can be arbitrary small, we get that $(k^{(s)})'(\hat \beta^{(s)})=r$ and therefore, $\beta^{(s)}\le \hat\beta^{(s)}$.

\hfill$\Box$

\skp\noi
{\bf Proof of Proposition \ref{prop_43}:} Fix $\eps\in(0,\iy)$. We start the proof by arguing uniqueness of the solution and then we move on to showing that a solution indeed exists.

{\bf Uniqueness:} We first argue that there cannot be two different $\calC^2([0,b],\R)$ solutions of HJB($\eps$). Notice that the non-linearity of our differential equation prevents us from using the same proof given in \cite[Proposition 2.2]{ata-shi}. Arguing by contradiction, assume that there are two $\calC^2([0,b],\R)$ solutions $f_1\ne f_2$. Without loss of generality, assume that there exists $x\in[0,b]$ such that $f_1(x)>f_2(x)$. Take $y_0\in\argmax_{x\in[0,b]}\{f_1(x)-f_2(x)\}$. Clearly, $f_1(y_0)>f_2(y_0)$. Also, $f_1'(0)=f_2'(0)=0$ and $f_1'(b)=f_2'(b)=r$ since they are both solutions to HJB($\eps$). Now, if $y_0\in(0,b)$, then first order condition for $f_1-f_2$ implies that $f_1'(y_0)=f_2'(y_0)$, and therefore it holds whether $y_0$ is an internal point or an endpoint of $[0,b]$.

By the structure of HJB($\eps$) applied to $f_1$, we get 
\begin{align}\label{418}
	f_1''(y_0)+H(y_0,f_1(y_0),f_1'(y_0))\ge 0.
\end{align}
The same inequality holds for $f_2$. Now, if 
\begin{align}\label{419}
	f_2''(y_0)+H(y_0,f_2(y_0),f_2'(y_0))=0,\end{align}
then by subtracting \eqref{419} form \eqref{418}, and recalling that $f_1(y_0)>f_2(y_0)$ and $f_1'(y_0)=f_2'(y_0)$, we get that $f_1''(y_0)>f_2''(y_0)$. Therefore, in a small neighborhood of $y_0$ (in case that $y_0=0$ or $y_0=b$, we take right- and left-neighborhoods, respectively), $f_1(x)-f_2(x)>f_1(y_0)-f_2(y_0)$, which contradicts the definition of $y_0$. 

Consider now the case
\begin{align}\notag
	f_2''(y_0)+H(y_0,f_2(y_0),f_2'(y_0))>0.
\end{align}
Since $f_2$ is a $\calC^2([0,b],\R)$ solution of HJB($\eps$), it follows that $0\le y_1<y_0<y_2\le b$, where
\begin{align}\notag
	y_1&:=\sup\left\{x\in[0,y_0] : f_2''(x)+H(x,f_2(x),f_2'(x))=0\right\}\vee 0,\\\notag
	y_2&:=\inf\left\{x\in[y_0,b] : f_2''(x)+H(x,f_2(x),f_2'(x))=0\right\}\wedge b.
\end{align}
We use the convention that $\sup\emptyset=-\inf\emptyset=-\iy$.  
Notice that on the interval $(y_1,y_2)$, which includes $y_0$, $f_2''(x)+H(x,f_2(x),f_2'(x))>0$. Since $f_2$ solves HJB($\eps$), it follows that $f_2'=0$ on $[y_1,y_2]$ or $f_2'=r$ on $[y_1,y_2]$. In the former case, we get $f_2(y_0)=f_2(y_2)$ and since $f_2'(b)=r$, $y_2<b$. Since $f_1'\ge 0$ as a solution of HJB($\eps$), we get that $f_1(y_2)-f_2(y_2)\ge f_1(y_0)-f_2(y_0)$, which by the definition of $y_0$ is in fact an equality. Therefore, $y_2 \in\argmax_{x\in[0,b]}\{f_1(x)-f_2(x)\}$. Since $y_2<b$ it follows by the smoothness of $f_2$ and the definition of $y_2$ that \eqref{418} and \eqref{419} hold with $y_2$ replacing $y_0$. Repeating the same argument proceeding \eqref{419}, with $y_2$ replacing $y_0$ yields a contradiction. Now if, $f_2'=r$ on $[y_1,y_2]$, a similar argument hold by using $f_1'\le r$ and the point $y_1$ instead of $y_2$.

\skp{\bf Existence:} 
We now construct a $\calC^2([0,b],\R)$ solution of \eqref{414} with some $\beta_\eps\in(0,b]$ that also solves HJB($\eps$). 
The solution is shown to satisfy also $k_{\beta_\eps}'\le r$. Repeating the same arguments given in \cite[Proposition 6.1]{har-tak}, we get that $k_{\beta_\eps}''(x)+H(x,k_{\beta_\eps}(x),k_{\beta_\eps}'(x))\ge 0$ on $[\beta,b]$ and therefore, $k_{\beta_\eps}$ satisfies HJB($\eps$) and we are done.

To solve \eqref{414}, we hereby consider the Cauchy problem given in \eqref{423}. Since eventually we find a solution to \eqref{423} that satisfies $0\le (k^{(s)}) '\le r$, it also solves the same ordinary differential equation from \eqref{423} with $H$ replacing $H_F$. The rest of the proof is performed in two steps. First, we prove the existence of $s^*\in\R$ for which the parameter $\beta^{(s^*)}\in(0,b]$  from \eqref{ac09} satisfies the following
\begin{align}\label{ac01}
	\text{$(k^{(s^*)})'(x)<r\;$ on $\;[0,\beta^{(s^*)})$},
\end{align}
\begin{align}\label{ac01b}
	\text{$\quad (k^{(s^*)})'(\beta^{(s^*)})=r$},
\end{align}
and in case that $\beta^{(s^*)}<b$, also 
\begin{align}\label{ac02}
	(k^{(s^*)})''(\beta^{(s^*)})=0.
\end{align} 
In the second step we show that $(k^{(s)})'\ge 0$ on $[0,\beta^{(s^*)}]$. Therefore, the function
\begin{align}\notag
	k(x)=\begin{cases}
		k^{(s^*)}(x),\qquad &0\le x< \beta^{(s^*)},\\\notag
		k^{(s^*)}(\beta^{(s^*)})+r(x-\beta^{(s^*)}),\qquad &\beta^{(s^*)}\le x\le b,
	\end{cases}
\end{align}
satisfies \eqref{414}, and the proof is done by setting 
\begin{align}\label{424a}
	\beta_\eps:=\beta^{(s^*)}\qquad\text{and}\qquad k_{\beta_\eps}:=k^{(s^*)}.
\end{align}
As a conclusion, we get that 
\begin{align}\label{424aa}
	s^*=k_{\beta_\eps}(0)=V(0;\eps).
\end{align}

{\bf Step 1:} 
Notice that \eqref{ac01} holds trivially for every $s^*$ by the definition of $\beta^{(s)}$. Therefore, we only need to check that \eqref{ac01b} and \eqref{ac02} hold. 

The first observation is that for sufficiently large $s$, $\beta^{(s)}<b$ and $(k^{(s)})''(\beta^{(s)})>0$. Also, for sufficiently small $s$, $\beta^{(s_0)}=b$ and $(k^{(s_0)})'(b)<r$. We prove the first part, the second one follows by the same lines and is therefore omitted. 
Fix $s_1> (2|m|r+2\sigma^2\eps r^2+h(b)+r\sigma^2/b)/\varrho$. Now, the function $k^{(s_1)}$ is nondecreasing on $[0,b]$. Indeed, arguing by contradiction, assume that there is $x\in(0,b)$ such that $(k^{(s_1)})'(x)<0$, then consider $y_3:=\inf\{x\in[0,b] : (k^{(s_1)})'(x)<0\}$. The smoothness of $k^{(s_1)}$ on $[0,b]$ and the initial condition $(k^{(s_1)})'(0)=0$ imply that $(k^{(s_1)})'(y_3)=0$. So, 
\begin{align}\label{426}
	\frac{1}{2}\sigma^2(k^{(s_1)})''(y_3)=\varrho k^{(s_1)}(y_3)-h(y_3)\ge \varrho k^{(s_1)}(0)-h(b)=\varrho s_1-h(b)>0.
\end{align}
The first inequality follows since $k^{(s_1)}$ is nondecreasing on $[0,y_3]$. This is a contradiction to the definition of $y_3$, since for every sufficiently small $\delta>0$ we have $(k^{(s_1)})'(y_3+\delta)>0$. 
Therefore, $(k^{(s_1)})'\ge 0$ on $[0,b]$. Combining it with \eqref{423}, $|F|\le 2r$, $k^{(s_1)}(0)=s_1$, the choice of $s_1$, and that $h$ is increasing (see the paragraph below \eqref{319b}), we get that for every $x\in[0,b]$,
\begin{align}\notag
	(k^{(s_1)})''(x)&=-\frac{2}{\sigma^2}\left(mF((k^{(s_1)})'(x))+\frac{1}{2}\sigma^2\eps F^2((k^{(s_1)})'(x))-\varrho k^{(s_1)}(x)+h(x)\right)\\\notag
	&> \frac{2}{\sigma^2}\left(\varrho k^{(s_1)}(0)-2|m|r-\frac{1}{2}\sigma^2r^2\eps-h(b)\right)=\frac{2r}{b}
	.
\end{align}
In particular $(k^{(s_1)})''(\beta^{(s_1)})>0$. Also $(k^{(s_1)})'(b/2)=\int_0^{b/2} (k^{(s_1)})''(u)du\ge r$ and therefore, $\beta^{(s_1)}\le b/2<b$. 

By the choices of $s_0$ and $s_1$, the following infimum is attained. Set 
\begin{align}\notag
	s^*:=&\inf\{s\in(-s_0,s_1) : \forall u\in(s,s_1),\; \beta^{(u)}<b\text{  and  }(k^{(u)})''(\beta^{(u)})>0\}.
\end{align}
In case that $\beta^{(s^*)}=b$, then by the definition of $s^*$, for every $u\in(s^*,s_1)$, $\beta^{(u)}<b$ and therefore, from \eqref{ac04} we have $(k^{(u)})'(\beta^{(u)})=r$. Now, the continuity of the mapping $s\mapsto (k^{(s)})'(\beta^{(s)})$ at $s=s^*$ in this case, provided in Lemma \ref{lem_42} implies that $(k^{(s^*)})'(\beta^{(s^*)})=r$ and \eqref{ac01b} holds in this case.

In case that $\beta^{(s^*)}<b$ then \eqref{ac04} implies that \eqref{ac01b} holds. We now claim that \eqref{ac02} holds, that is, $(k^{(s^*)})''(\beta^{(s^*)})=0$. 
%
%
Otherwise, if $(k^{(s^*)})''(\beta^{(s^*)})<0$, then by the continuity of the mapping $s\mapsto (\beta^{(s)},(k^{(s)})''(\beta^{(s)}))$ in this case (see Lemma \ref{lem_42}), we get that for every sufficiently small $\nu>0$, $\beta^{(s^*+\nu)}<b$ and $(k^{(s^*+\nu)})''(\beta^{(s^*+\nu)})<0$, which contradicts the definition of $s^*$. If however, $(k^{(s^*)})''(\beta^{(s^*)})>0$ then Lemma \ref{lem_42} again implies that for every sufficiently small $\nu>0$, one has $\beta^{(s^*-\nu)}<b$ and $(k^{(s^*-\nu)})''(\beta^{(s^*-\nu)})>0$, in contradiction to the definition of $s^*$. Therefore, $(k^{(s^*)})''(\beta^{(s^*)})=0$ and \eqref{ac02} holds. 

{\bf Step 2:} In this step we show that 
\begin{align}\label{428}
	(k^{(s^*)})'(x)\ge 0,\qquad x\in[0,\beta^{(s^*)}]. 
\end{align}
Arguing by contradiction, assume that there is $y_5\in(0,\beta^{(s^*)})$ with $(k^{(s^*)})'(y_5)< 0$. We omitted the endpoints where the derivatives are $0$ and $r$. The smoothness of the function $k^{(s^*)}$ together with $(k^{(s^*)})'(0)=0$ and $(k^{(s^*)})'(\beta^{(s^*)})=r$ (from the previous step) imply that the following supremum and infimum are attained 
\begin{align}\notag
	y_4&:=\inf\{x\in[0,y_5) : \forall y\in(x,y_5), \; (k^{(s^*)})'(y)\le0\},\\\notag
	y_6&:=\sup\{x\in(y_5,\beta^{(s^*)}) : \forall y\in(x,y_5), \; (k^{(s^*)})'(y)\le0\}.
\end{align}
Also, $(k^{(s^*)})'(y_4)=(k^{(s^*)})'(y_6)=0$ and
$(k^{(s^*)})''(y_4)\le0\le (k^{(s^*)})''(y_6)$. Substituting these relations in \eqref{423}, we obtain,
\begin{align}\label{430}
	\varrho k^{(s^*)}(y_4)-h(y_4)\le 0\le \varrho k^{(s^*)}(y_6)-h(y_6).
\end{align}
Since $(k^{(s^*)})'(x)\le0$ on $[y_4,y_6]$, and by the smoothness of $k^{(s^*)}$, there is a subinterval on which $(k^{(s^*)})'(x)<0$, one has $k^{(s^*)}(y_4)> k^{(s^*)}(y_6)$. Recall that the function $h$ is increasing (see the paragraph below \eqref{319b}) and therefore, $h(y_4)< h(y_6)$. The last two inequalities contradict \eqref{430} and therefore, \eqref{428} holds.

\hfill$\Box$

%
%

\section{Optimal strategies and equilibria in the games}\label{sec_5}
\beginsec

In Theorem \ref{thm_41} we claimed that the minimizer has an optimal strategy, which is a $\beta_\eps$-reflecting strategy with $\beta_\eps$ given in \eqref{440}. However, we did not argue uniqueness in the sense that for every $\beta\ne\beta_\eps$, the $\beta$-reflecting strategy is strictly sub-optimal. Neither did Harrison and Taksar  \cite{har-tak} nor Atar and Shifrin \cite{ata-shi}. Hence, in the following discussion, which evolves around this issue, we allow $\eps=0$ and consider $\eps\in[0,\iy)$. 

Recall the definition of $\beta_\eps$ from \eqref{440}. Define now,
\begin{align}\label{441}
	\hat\beta_\eps&:=\sup\left\{x\in(0,b] : \forall y\le x,\;
	V''(y;\eps)+H(y,V(y;\eps),V'(y;\eps))=0\right\}.
\end{align}
From \eqref{414z} it follows that $\hat\beta_\eps\ge\beta_\eps$. Since on the interval $[\beta_\eps,\hat\beta_\eps]$ both of the conditions 
\begin{align}\label{442}
	V''(x;\eps)+H(x,V(y;\eps),V'(x;\eps))=0\qquad
	\text{and}\qquad
	V'(x;\eps)=r
\end{align} hold, it follows that $V$ solves \eqref{414} for any $\beta\in[\beta_\eps,\hat\beta_\eps]$. Proposition \ref{prop_42} implies that every such $\beta$-reflecting strategy is optimal. 
Thus, the non-uniqueness of the optimal reflecting strategy is equivalent to the existence of non-degenerate interval $[\beta_\eps,\hat\beta_\eps]$, on which the equations in \eqref{441} hold. Combining both of them, we get that $V(x;\eps)=(mr+\sigma^2r^2\eps/2+h(x))/\varrho$. Using again $V'(x;\eps)=r$, we get that $h'(x)=\varrho r$. Recall that $h$ is piecewise linear with slopes in the set $\{h_1\mu_1,\ldots,h_I\mu_I\}$ (see the paragraph below \eqref{319b}). Thus, if $\varrho r$ does not belong to this set, we get uniqueness of the optimal reflecting strategy. Another sufficient condition would be that $m$ is sufficiently negative such that $(mr+\sigma^2r^2\eps/2+h(b))/\varrho\le 0$, since the value function cannot be non-positive. The arguments above are summarized in the following proposition.
\begin{proposition}\label{prop_44}
	Fix $\eps\in[0,\iy)$. A $\beta$-reflecting strategy is optimal if and only if $\beta\in[\beta_\eps, \hat\beta_{\eps}]$. If 
	\begin{align}\label{ac10}
		\text{for every $i\in[I]$, $\hat h_i\mu_i\neq \varrho r\quad$ or $\quad mr+\frac{1}{2}\sigma^2\eps r^2+h(b)\le 0$}, 
	\end{align}
	then $\beta_\eps= \hat\beta_{\eps}$ and there is a unique optimal reflecting strategy.
\end{proposition}
Although there might be weaker assumptions under which uniqueness hold, we do not aim to find them, and remain contented with the not so restrictive conditions given in \eqref{ac10}. 

We already know that the minimizer has an optimal reflecting strategy. We now show that for every $\kaboom\in(0,\iy)^I$ and $\eps$ given through \eqref{310g}, both the RSDG and the MSDG have equilibria. Specifically, each player has an optimal strategy that is good against any strategy of the other player. 
Recall that Proposition \ref{prop_31} connects between the costs in the games. Fix $\hat x_0\in\calX$, $\kaboom\in(0,\iy)^I$, $\bar\beta_\eps\in[\beta_\eps,\hat\beta_\eps]$, and an $I$-dimensional standard Brownian motion $\hat B$. Set $B=\sigma^{-1}\theta\hat\sigma\cdot\hat B$ and recall \eqref{310b} and \eqref{310g}. 
Let $(Y_{\bar\beta_\eps},R_{\bar\beta_\eps})$ be an optimal reflecting strategy for the minimizer. Also, let $\Q_V=\Q_{V(\cdot;\eps)}$ be the measure driven by $\psi_V=\psi_{V(\cdot;\eps)}$ (see \eqref{407}). Recall the definition of the function $\gamma$ from \eqref{323} and the relations given in \eqref{325}--\eqref{326}. 
Set the strategy $(\hat Y_{\bar\beta_\eps},\hat R_{\bar\beta_\eps})$ by
\begin{align}\notag
	\hat R_{\bar\beta_\eps}(t):=R_{\bar\beta_\eps}(t)\mu_{i^*}e_{i^*}\qquad\text{and}\qquad
	\hat Y_{\bar\beta_\eps}(t):=\hat X_{\bar\beta_\eps}(t)-\hat x_0-\hat mt-\hat\sigma \hat B(t) +\hat R_{\bar\beta_\eps}(t),\quad t\in\R_+,
\end{align}
where $i^*$ is given in \eqref{320} and for any $t\in\R_+$,
\begin{align}\notag
	\hat X_{\bar\beta_\eps}(t):&=\gamma (X_{\bar\beta_\eps}(t)),\\\notag
	X_{\bar\beta_\eps}(t)&= x_0+ mt+\sigma  B(t)+ Y_{\bar\beta_\eps}(t)- R_{\bar\beta_\eps}(t).
\end{align}
Moreover, let $\hat\Q_V$ be the measure associated with $\hat\psi_V(t)=(\hat\psi_{V,1}(t),\ldots,\hat\psi_{V,I}(t))$, given by
\begin{align}\label{444d}
	\hat\psi_{V,i}(t):=\frac{\sigma\psi_V(t)(\theta\hat\sigma)_i\kaboom_i}{\sum_{j=1}^I(\theta\hat\sigma)^2_j\kaboom_j}
	,\qquad i\in[I],\quad t\in\R_+.
\end{align}

We now claim that the mentioned strategies form equilibria in the games. Using the $\bar \beta_\eps$-reflecting strategy in the RSDG it follows that the associated workload dynamics $X_{\bar \beta_\eps}$ is a reflected diffusion that moves continuously along the interval $[0,\bar \beta_\eps]$. The function $\gamma$ continuously maps the workload process to the bold curve from Figure \ref{fig_1}. If for example $m=\bar \beta_\eps$ in figure \ref{fig_1}, then the process $\hat X_{\beta_\eps}$ moves continuously along the bold curve between the points $J$ and $M$, where it is reflected at these two points.

The next theorem generalizes Corollary  \ref{cor_31} and Theorem \ref{thm_41} and provides equilibria in the games.
\begin{theorem}\label{thm_42}
	Fix $\eps\in(0,\iy)$. Using the notation above, the triplets $(Y_{\bar\beta_\eps},R_{\bar\beta_\eps},\Q_V)$ and $(\hat Y_{\bar\beta_\eps},\hat R_{\bar\beta_\eps},\hat\Q_V)$ form equilibria in the RSDG and the MSDG, respectively. That is,
	\begin{align}\notag
		&V(x_0;\eps)=
		\sup_{\Q\in\calQ(x_0)}J(x_0,Y_{\bar\beta_\eps},R_{\bar\beta_\eps},\Q;\eps)=
		J(x_0,Y_{\bar\beta_\eps},R_{\bar\beta_\eps},\Q_V;\eps)=
		\inf_{(Y,R)\in\calA(x_0)}J(x_0,Y,R,\Q_V;\eps)\\\notag
		&\quad =\hat V(\hat x_0;\kaboom)=
		\sup_{\hat \Q\in\hat \calQ(\hat x_0)}\hat J(x_0,\hat Y_{\bar\beta_\eps},\hat R_{\bar\beta_\eps},\hat \Q;\kaboom)=
		\hat J(\hat x_0,\hat Y_{\bar\beta_\eps},\hat R_{\bar\beta_\eps},\hat \Q_V;\kaboom)=
		\inf_{(\hat Y,\hat R)\in\hat \calA(\hat x_0)}\hat J(\hat x_0,\hat Y,\hat R,\hat \Q_V;\kaboom).
	\end{align}
\end{theorem}

\begin{remark}\label{rem_ac01}
	Notice that the optimal policies for the minimizer in the MSDG has the same structure as the optimal policy in the BCP given in \cite{ata-shi}. The only difference emerge from the cut-off point $\beta_\eps$, which affects the point of reflection on the curve $\gamma$. Such a result is not  obvious due to the non-stationarity structure of the problem caused by the existence of the maximizer player. 
	Furthermore, the structures of the equilibria share similarities with the optimal policies in the differential game given in \cite{ata-coh}. More specifically, \cite[Theorem 3.3]{ata-coh} states that the minimizer's optimal strategy in the one-dimensional deterministic differential game is a reflecting strategy (called there a `barrier strategy'). Under optimality, the maximal player in the same game uses a drift change that is illustrated 
	more compactly in \cite[Section 3, (30)--(31)]{ata-coh2017}. The relationship between the multidimensional and the one-dimensional games is given in \cite[Appendix A]{ata-coh}.
\end{remark}

\noi{\bf Proof of Theorem \ref{thm_42}:} 
In this proof we do not check admissibility. It can be verified easily the same way as in Proposition \ref{prop_31} using Corollary \ref{cor_42b}. The first two equalities follow from the optimality of $(Y_{\bar\beta_\eps},R_{\bar\beta_\eps})$ and from \eqref{415}. The third equality follows since on the one hand, by the definition of $V$,
\begin{align}\notag
	V(x_0;\eps)\ge\inf_{(Y,R)\in\calA(x_0)}J(x_0,Y,R,\Q_V;\eps)
\end{align}
and on the other hand, the reversed inequality follows from Proposition \ref{prop_41} since by Theorem \ref{thm_41}, $V\in\calC^2([0,b],\R)$.
The forth equality follows by Corollary \ref{cor_31}. For the fifth and sixth equalities notice that $\hat V(\hat x_0;\kaboom)=V(x_0;\eps)$ and that 
\begin{align}\notag
	\hat V(\hat x_0;\kaboom)\le \sup_{\hat\Q\in\hat\calQ(\hat x_0)}\hat J(\hat x_0,\hat Y_{\bar\beta_\eps},\hat R_{\bar\beta_\eps},\hat\Q;\kaboom)\le \sup_{\Q\in\calQ(x_0)}J(x_0,Y_{\bar\beta_\eps},R_{\bar\beta_\eps},\Q;\eps)=V(x_0;\eps),
\end{align}
where the second inequality follows by \eqref{327c}. 
For the last equality, notice that we already established
\begin{align}\notag
	\hat V(\hat x_0;\kaboom)=\sup_{\hat \Q\in\hat \calQ(\hat x_0)}\hat J(x_0,\hat Y_{\bar\beta_\eps},\hat R_{\bar\beta_\eps},\hat \Q;\kaboom)\ge
	\hat J(\hat x_0,\hat Y_{\bar\beta_\eps},\hat R_{\bar\beta_\eps},\hat \Q_V;\kaboom)\ge
	\inf_{(\hat Y,\hat R)\in\hat \calA(\hat x_0)}\hat J(\hat x_0,\hat Y,\hat R,\hat \Q_V;\kaboom).
\end{align}
Therefore, it is sufficient to show that
$\inf_{(\hat Y,\hat R)\in\hat \calA(\hat x_0)}\hat J(\hat x_0,\hat Y,\hat R,\hat \Q_V;\kaboom)\ge V(x_0;\eps)$, which holds by \eqref{324a} (replace the subindex $*$ with $V$) as follows,
\begin{align}\notag
	\inf_{(\hat Y,\hat R)\in\hat \calA(\hat x_0)}\hat J(\hat x_0,\hat Y,\hat R,\hat\Q_V;\kaboom)\ge \inf_{( Y, R)\in \calA( x_0)}J(x_0,Y,R,\Q_V;\eps)=V(x_0;\eps).
\end{align}
%
%

\hfill$\Box$

\section{Dependency on the ambiguity parameters}\label{sec_6}
\beginsec

In this section we study the dependence of the value functions and the optimal cut-offs on the ambiguity parameters. 
We show continuity and that as $\eps\to0$ our model converges to the risk-neutral model, studied by Harrison and Taksar \cite{har-tak} and by Atar and Shifrin \cite{ata-shi}. For this, recall the definition of $V(\cdot;0)$ given in Remark \ref{rem_32}.

\begin{theorem}\label{thm_45}
	The mapping $[0,\iy)\ni\eps\mapsto (V(\cdot;\eps),V'(\cdot;\eps))$ is increasing and continuous in the uniform norm topology taken on the interval $[0,b]$. 
	Moreover,  there is a constant $C>0$ such that for every $\eps\in(0,\iy)$, $\sup_{x\in[0,b]}|V(x;\eps)-V(x;0)|\le C\eps$. Also, $\lim_{\eps\to\iy}V(\cdot;\eps)=\iy$, uniformly on $[0,b]$. 
	Finally, consider the relations $\prec$ and $\preceq$ on $(0,\iy)^I\times(0,\iy)^I$ given by,
	\begin{align}\notag
		\kaboom\prec\kaboom'\quad\text{(resp., $\preceq$)}\qquad\text{if and only if}\qquad \sum_{i=1}^I(\theta\hat\sigma)^2_i\kaboom_i< \sum_{i=1}^I(\theta\hat\sigma)^2_i\kaboom_i'\quad\text{(resp., $\le$)}.
	\end{align}
	Then the mapping $(0,\iy)^I\ni\kaboom\mapsto\hat V(\cdot;\kaboom)$ is increasing w.r.t.~$\prec$ and continuous in the uniform norm topology taken on $\calX$.
\end{theorem}

\noi{\bf Proof:} The last part of the theorem merely follows by the first one, \eqref{310g}, and Corollary \ref{cor_31}. Therefore it is omitted and we turn to proving the first part of the theorem. We start by showing the monotonicity and continuity of the mapping $(0,\iy)\ni\eps\mapsto V(\cdot;\eps)$. The proof for $\eps=0$ is given separately. 
Fix $0<\eps_2<\eps_1$. Denote by $\Q_i=\Q_{V(\cdot;\eps_i)}$ and $\psi_i=\psi_{V(\cdot;\eps_i)}$, $i=1,2$. Then, for every $x_0\in[0,b]$ one has
\begin{align}\notag
	V(x_0;\eps_1)&=\inf_{(Y,R)\in\calA(x_0)}J(x_0,Y,R,\Q_1;\eps_1) \\\notag
	&= \inf_{(Y,R)\in\calA(x_0)}\left[J(x_0,Y,R,\Q_1;\eps_2)+\frac{1}{2}\left(\frac{1}{\eps_2}-\frac{1}{\eps_1}\right) \int_0^\iy e^{-\varrho t}\E^{\Q_1}\left[\psi_1^2(t)\right]dt\right]\\\notag
	&\le \inf_{(Y,R)\in\calA(x_0)}\left[J(x_0,Y,R,\Q_1;\eps_2)+\frac{\eps_1\sigma^2r^2}{2\eps_2\varrho}(\eps_1-\eps_2)\right]\\\notag
	&\le \inf_{(Y,R)\in\calA(x_0)}\sup_{\Q\in\calQ(x_0)}\left[J(x_0,Y,R,\Q;\eps_2)\right]+\frac{\eps_1\sigma^2r^2}{2\eps_2\varrho}(\eps_1-\eps_2)\\\notag
	&=V(x_0;\eps_2)+\frac{\eps_1\sigma^2r^2}{2\eps_2\varrho}(\eps_1-\eps_2).
\end{align}
The first equality follows by Theorem \ref{thm_42}. The second equality follows by  \eqref{321}. The first inequality follows since by \eqref{407}, $\psi_1(t)=\eps\sigma V'(X(t);\eps_1)$ and since $V'(x;\eps_1)\le r$, see Corollary \ref{cor_42b}. The second inequality is trivial and finally, the last equality follows by the definition of $V$, see \eqref{322}.

On the other hand,
\begin{align}\notag
	V(x_0;\eps_1)&\ge\inf_{(Y,R)\in\calA(x_0)}J(x_0,Y,R,\Q_2;\eps_1) \\\notag
	&= \inf_{(Y,R)\in\calA(x_0)}\left[J(x_0,Y,R,\Q_2;\eps_2)+\frac{1}{2}\left(\frac{1}{\eps_2}-\frac{1}{\eps_1}\right) \int_0^\iy e^{-\varrho t}\E^{\Q_2}\left[\psi_2^2(t)\right]dt\right]\\\notag
	&> \inf_{(Y,R)\in\calA(x_0)}J(x_0,Y,R,\Q_2;\eps_2)\\\notag
	&=V(x_0;\eps_2).
\end{align}
The first inequality follows by the definition of $V$. The first equality follows by \eqref{321}. The strict inequality follows since $\eps_1>\eps_2$ and since by \eqref{407}, $\psi_1(t)=\eps\sigma V'(X(t);\eps_1)$ and since $V'(x;\eps_1)\ge 0$, see Corollary \ref{cor_42b}. Obviously, with $\Q_2$-probability zero, $V'(X(t);\eps_1)=0$ for almost every $t\in\R_+$ w.r.t.~Lebesgue measure. The last equality follows by Theorem \ref{thm_42}.
Combining the last two sets of relations, we obtain,
\begin{align}\notag
	V(x_0;\eps_2)< V(x_0,;\eps_1)\le V(x_0;\eps_2)+\frac{\eps_1\sigma^2r^2}{2\eps_2\varrho}(\eps_1-\eps_2)
\end{align}
and the monotonicity and continuity of $\eps\mapsto V(\cdot;\eps)$ is proven on the interval $(0,\iy)$. The monotonicity at $\eps=0$ follows by the following sequence of relations,
\begin{align}\label{453}
	V(x_0;0)&=\inf_{(Y,R)\in\calA(x_0)}\E^{ \PP}\Big[\int_0^\iy e^{-\varrho t}( h(  X(t))dt + rd R(t)) \Big]\\\notag
	&=\inf_{(Y,R)\in\calA(x_0)}J(x_0,Y,R,\PP;\eps_1)\\\notag
	&\le\inf_{(Y,R)\in\calA(x_0)}\sup_{\Q\in\calQ(x_0)}J(x_0,Y,R,\Q;\eps_1)\\\notag
	&=V(x_0;\eps_1),
\end{align}
where $x_0\in[0,b]$ and $\eps_1\in(0,\iy)$. The first equality follows by the definition of the cost in the risk-neutral case, see \eqref{no_amb}. The second equality follows since $L^\varrho(\PP\|\PP)=0$.

We now turn to proving the continuity at $\eps=0$. Notice that in the arguments above, the inequality $\eps_2>0$ cannot be relaxed to $\eps_2\ge 0$ since we divide by $\eps_2$. Therefore, we come up with another proof for continuity at $0$. 
Recall that the arguments in Section \ref{sec_5} included the case where $\eps$ is zero. Fix an optimal cut-off in the risk-neutral problem $\bar\beta_0\in[\beta_0,\hat\beta_0]$. Fix also $\eps>0$ and set $\Q^\eps:=\Q_{V(\cdot;\eps)}$ and $\psi^\eps:=\psi_{V(\cdot;\eps)}$. 
From \eqref{453} and Theorem \ref{thm_42} one has,
\begin{align}\notag
	V(x_0;0)\le V(x_0;\eps)= J(x_0,Y_{\bar\beta_0},R_{\bar\beta_0},\Q^\eps;\eps).
\end{align}
Recall that by \eqref{321} and \eqref{407},
\begin{align}\label{453a}
	&J(x_0,Y_{\bar\beta_0},R_{\bar\beta_0},\Q^\eps;\eps)\\\notag
	&\quad=\E^{ \Q^\eps}\Big[\int_0^\iy e^{-\varrho t}( h(  X_{\bar\beta_0}(t))dt + rd R_{\bar\beta_0}(t)) \Big]-\eps\sigma^2\int_0^\iy e^{-\varrho t}(V'(X_{\bar\beta_0}(s);\eps))^2ds,
\end{align}
where 
\begin{align}\label{453b}
	X_{\bar\beta_0}(t)= x_0+ mt+\sigma  B(t)+ Y_{\bar\beta_0}(t)- R_{\bar\beta_0}(t),\qquad t\in\R_+.
\end{align}
Since $0\le V'(\cdot;\eps)\le r$
, the last term on the r.h.s.~of \eqref{453a} goes to zero as $\eps\to 0$. Hence, in order to show the limit $\lim_{\eps\to 0}V(x_0;\eps)=V(x_0;0)$ it is sufficient to show that
\begin{align}\label{455}
	\lim_{\eps\to 0}\;\E^{ \Q^\eps}\Big[\int_0^\iy e^{-\varrho t}( h(  X_{\bar\beta_0}(t))dt + rd R_{\bar\beta_0}(t)) \Big]=\E^{ \PP}\Big[\int_0^\iy e^{-\varrho t}( h(  X_{\bar\beta_0}(t))dt + rd R_{\bar\beta_0}(t)) \Big].
\end{align}

The proof of \eqref{455} relies on the continuity property of the Skorokhod mapping in addition to a coupling argument. We consider a probability space $(\Omega,\calF,\{\calF_t\},\check\PP)$ that supports a one-dimensional standard Brownian motion $B$, adapted to the filtration $\{\calF_t\}$, such that for every $0\le s<t$, $B(t)-B(s)$ is independent of $\calF_s$ under $\check \PP$. Recall Definition \ref{def_Skorokhod} and consider the following processes,
\begin{align}\notag
	(X^\eps,Y^\eps,R^\eps)(t)&:=\Gamma_{[0,\bar\beta_0]}\left(x_0+m\cdot+\int_0^\cdot\eps \sigma^2V'(X^\eps(s);\eps)(s)ds+\sigma B(\cdot)\right)(t),\\\notag
	(X^0,Y^0,R^0)(t)&:=\Gamma_{[0,\bar\beta_0]}\left(x_0+m\cdot+\sigma B(\cdot)\right)(t).
\end{align}
Notice that by \eqref{407}, $(X^0,Y^0,R^0)$ (resp., $(X^\eps,Y^\eps,R^\eps)$) has the same distribution under the measure $\check\PP$ as $(X_{\bar\beta_0},Y_{\bar\beta_0},R_{\bar\beta_0})$ given in \eqref{453b} under the measure $\PP$ (resp., $\Q^\eps$, see \eqref{317}). Hence, \eqref{455} follows once we show that 
\begin{align}\label{457}
	\lim_{\eps\to 0}\;\E^{ \check\PP}\Big[\int_0^\iy e^{-\varrho t}( h(  X^\eps(t))dt + rd R^\eps(t)) \Big]=\E^{ \check\PP}\Big[\int_0^\iy e^{-\varrho t}( h(  X^0(t))dt + rd R^0(t)) \Big].
\end{align}
Now,
\begin{align}\label{458}
	&\left|\E^{ \check\PP}\Big[\int_0^\iy e^{-\varrho t}( h(  X^\eps(t))dt + rd R^\eps(t)) \Big]-\E^{ \check\PP}\Big[\int_0^\iy e^{-\varrho t}( h(  X^0(t))dt + rd R^0(t)) \Big]\right|\\\notag
	&\quad=\left|\E^{ \check\PP}\Big[\int_0^\iy e^{-\varrho t}[ h(  X^\eps(t))-h(X^0(t))]dt + \left.\int_0^\iy \varrho re^{-\varrho t}[R^\eps(t)-R^0(t)]dt +[e^{-\varrho t}(R^\eps(t)-R^0(t))]\right|_0^\iy \Big]\right|\\\notag
	&\quad\le
	\E^{ \check\PP}\Big[\int_0^\iy e^{-\varrho t}| h(  X^\eps(t))-h(X^0(t))|dt + \int_0^\iy \varrho re^{-\varrho t}|R^\eps(t)-R^0(t)|dt \Big].
\end{align}
The last inequality follows since $\lim_{t\to\iy} e^{-\varrho t}R^\eps(t)=0$, $\check\PP$-a.s.~and similarly for $R^0$. Indeed, consider $w(\cdot)=x_0+m\cdot+\int_0^\cdot\eps \sigma^2V'(X^\eps(s);\eps)(s)ds+\sigma B(\cdot)$ and $\tilde w=0$ in Lemma \ref{lem_Skorokhod}. Recalling the bound $0\le V'(\cdot,\eps)\le r$, it follows that $R^\eps(t)\le c_S (x_0+(m+\eps\sigma^2 r)t+\sigma \sup_{0\le s\le t}|B(s)|)$ and the limit follows since $\lim_{t\to\iy} e^{-\varrho t} \sup_{0\le s\le t}|B(s)|=0$, $\check\PP$-a.s. Using the Lipschitz continuity of $h$, again $0\le V'(\cdot,\eps)\le r$, and Lemma \ref{lem_Skorokhod} with $w(\cdot)=x_0+m\cdot+\int_0^\cdot\eps \sigma^2V'(X^\eps(s);\eps)(s)ds+\sigma B(\cdot)$ and $\tilde w(\cdot)=x_0+m\cdot+\sigma B(\cdot)$, we get that for every $t\in\R_+$,
\begin{align}\notag
	\sup_{s\in[0,t]}\|(X^\eps,Y^\eps,R^\eps)(s)-(X^0,Y^0,R^0)(s)\|\le \eps c_S\sigma^2rt.
\end{align}
Recalling that $h$ is Lipschitz, we get that the r.h.s.~of \eqref{458} is bounded above by $C\eps(1/\varrho+r)/\varrho$,where the constant $C>0$ depends on $c_S$, $\sigma^2r$, and the Lipschitz constant of $h$.
This implies \eqref{457}. Recalling \eqref{453a} and the estimations above, we get the bound 
$$\sup_{x\in[0,b]}|V(x;\eps)-V(x;0|\le (\sigma^2r^2+C\eps(1/\varrho+r)/\varrho)\eps$$ as required.

We now turn to proving the limit $\lim_{\eps\to\iy}V(\cdot;\eps)=\iy$. For this, consider the strategy of the maximizer given by $\Q_\eps$, which is associated with $\psi_\eps(s)=\eps^{1/4}$ for every $s\in\R_+$. Given that and recalling the presentation of the dynamics provided in \eqref{317}, it follows that, 
\begin{align}\label{ac11}
	X(t)= x_0+ (m+\sigma\eps^{1/4})t+ \sigma  B^{\Q_\eps}(t)+ Y(t)- R(t),\quad t\in\R_+,
\end{align}
where $ B^{ \Q_\eps}$ is an $\{\calF_t\}$-one-dimensional standard Brownian motion under $\Q_\eps$. His cost function is given by
\begin{align}\notag
	J( x_0, Y, R, \Q_\eps;\eps)&=
	\E^{ \Q}\Big[\int_0^\iy e^{-\varrho t}\Big( h( X(t))dt + r  d R(t)-\frac{1}{2\eps} \psi^2(t)dt\Big) \Big]\\\notag
	&=\E^{ \Q_\eps}\Big[\int_0^\iy e^{-\varrho t}( h( X(t))dt + r  d R(t))\Big]-\frac{1}{2\varrho^2\eps^{1/2}}.
\end{align}
Since the last term is constant, we can think as if the minimizer is facing a risk-neutral problem with the drift $mt$ replaced by $(m+\sigma\eps^{1/4})t$, where her  cost function has an additional deterministic term, which vanishes as $\eps\to\iy$. Therefore, we analyze only the limit of the expectation in the second line of the above. Consider an optimal $\beta$-reflecting strategy $(Y_\beta,R_\beta)$, which exists in the risk-neutral problem. From the dynamics in \eqref{ac11}, we get for every $t\in\R_+$,
\begin{align}\notag
	R_\beta(t)\ge x_0-X(t)+ (m+\sigma\eps^{1/4})t+ \sigma  B^{\Q_\eps}(t)+ Y_\beta(t)\ge -b+ (m+\sigma\eps^{1/4})t+ \sigma  B^{\Q_\eps}(t). 
\end{align}
The same arguments that lead to \eqref{458}, imply that
\begin{align}\notag
	&\E^{ \Q_\eps}\Big[\int_0^\iy e^{-\varrho t}( h( X(t))dt + r  d R_\beta(t))\Big]=
	\E^{ \Q_\eps}\Big[\int_0^\iy e^{-\varrho t}( h( X(t))+\varrho r R_\beta(t))dt\Big]\\\notag
	&\quad\ge\E^{ \Q_\eps}\Big[\int_0^\iy e^{-\varrho t}( h( X(t))+\varrho r (-b+ (m+\sigma\eps^{1/4})t+ \sigma  B^{\Q_\eps}(t)))dt\Big].
\end{align}
Standard estimates give that the last term goes to $\iy$ together with $\eps$.

\hfill$\Box$

\skp
From Theorem \ref{thm_42}, Corollary \ref{cor_42b}, and \eqref{407} it is easy to see the continuity of the maximizer's optimal control w.r.t.~$\eps$. Also, it is clear that as $\eps\to0$, $\psi_{V(\cdot;\eps)}\to0$, uniformly on $[0,b]$ and therefore, the maximizer's optimal strategy is relatively close to keep the original measure as is. Since we do not have explicit expressions for the cut-off points of the reflecting strategies $\beta_\eps$ and $\hat\beta_\eps$ described in Section \ref{sec_5}, the situation with the minimizer's optimal strategy is more subtle and is studied now. The next theorem provides a continuity property of the optimal cut-off points $\beta_\eps$ and $\hat\beta_\eps$ w.r.t.~$\eps$.
Recall the sufficient conditions given in \eqref{ac10} for uniqueness of the reflecting strategy.

\begin{theorem}\label{prop_46}
	For any given $\eps\in[0,\iy)$, 
	\begin{align}\label{460}
		\beta_\eps\le \liminf_{\delta\to 0}\beta_{\eps+\delta}\le \limsup_{\delta\to 0}\hat\beta_{\eps+\delta}&\le \hat\beta_\eps.
	\end{align}
	Hence, in case that $\hat\beta_\eps=\beta_\eps$ (for the given $\eps$), that is, if there is a unique optimal reflecting strategy, one has,
	\begin{align}\label{461}
		\lim_{\delta\to 0}\beta_{\eps+\delta}=\lim_{\delta\to 0}\hat\beta_{\eps+\delta}=\hat\beta_\eps.
	\end{align} 
\end{theorem}
As a preparation for the proof, we present an auxiliary function defined for every $\eps\in[0,\iy)$. Recall the definition of $k^{(s)}$ given in \eqref{423} and \eqref{424aa}. Set $l^{(\eps)}(x)=k^{(V(0;\eps))}(x)$, $x\in[0,b]$. More explicitly,
\begin{align}\label{462}
	\begin{cases}
		(l^{(\eps)})''(x)+H_F(x,l^{(\eps)}(x),(l^{(\eps)})'(x))=0,\qquad x\in[0,b],\\
		(l^{(\eps)})'(0)=0,\quad l^{(\eps)}(0)=V(0;\eps),
	\end{cases}
\end{align}
with the same $F$ given after \eqref{423}. 
The next lemma 
provides a continuity property of $l^{(\eps)}$ as a function of $\eps$.

\begin{lemma}\label{lem_43}
	The mapping $[0,\iy)\ni\eps\mapsto(l^{(\eps)}(\cdot),(l^{(\eps)})'(\cdot))$ is continuous in the uniform norm topology taken on the interval $[0,b]$.
\end{lemma}
{\bf Proof
	:} Fix $\eps_1,\eps_2\in\R_+$. Set $f_1=l^{(\eps_1)}$ and $f_2=l^{(\eps_2)}$. Then,
\begin{align}\notag
	f_1'(x)&=-\int_0^xH^{\eps_1}_F(y,f_1(y),f_1'(y))dy,\\\notag
	f_2'(x)&=-\int_0^xH^{\eps_2}_F(y,f_2(y),f_2'(y))dy = -\int_0^xH^{\eps_1}_F(y,f_2(y),f_2'(y))dy +\int_0^x(\eps_2-\eps_1)F(f_2'(y)),
\end{align}
where the index $\eps_i$ in $H^{\eps_i}_F$ emphasizes its dependence on $\eps_i$, $i=1,2$. 
From \eqref{ac05} it follows that there exists a constant $L>0$ independent of $\eps_1,\eps_2$ such that 
\begin{align}\notag
	|f_1'(x)-f_2'(x)|\le L\left(\int_0^x\Big[|f_1(y)-f_2(y)|+|f_1'(y)-f_2'(y)|\Big]dy +|\eps_1-\eps_2|\right).
\end{align}
Also, since $f_i(0)=V(0;\eps_i)$,
\begin{align}\notag
	|f_1(x)-f_2(x)|\le |V(0;\eps_1)-V(0;\eps_2)|+\int_0^x|f_1'(y)-f_2'(y)|dy.
\end{align}
From the last two inequalities and Gr\"{o}nwall's inequality we get that there is a constant $C>0$ independent of $\eps_1,\eps_2$, and $x$, such that for every $x\in[0,b]$
\begin{align}\notag
	|f_1(x)-f_2(x)|+|f_1'(x)-f_2'(x)|\le C\left( |V(0;\eps_1)-V(0;\eps_2)|+|\eps_1-\eps_2|\right).
\end{align}
Recalling from Theorem \ref{thm_45} that $\eps\mapsto V(\cdot;\eps)$ is continuous, we get that $\\$$\sup_{0\le x\le b}\left(|f_1(x)-f_2(x)|+|f_1'(x)-f_2'(x)|\right)$ is uniformly bounded by a function of $(\eps_1,\eps_2)$ that goes to zero as $(\eps_1-\eps_2)\to 0$. 

\hfill$\Box$

\skp\noi{\bf Proof of Theorem \ref{prop_46}:} 
For any $\eps'\in[0,\iy)$, the function $V(\cdot;\eps')$ satisfies \eqref{414z} and therefore also \eqref{462} on $x\in[0,\hat\beta_{\eps'}]$. Uniqueness of the solution implies that 
\begin{align}\label{467}
	V(x;\eps')=l^{(\eps')}(x), \qquad x\in[0,\hat\beta_{\eps'}], 
\end{align}
see \cite[Section 0.3.1]{Polyanin2003}. 

To get the first inequality in \eqref{460} notice that from \eqref{467}, the definition of $\beta_\eps$ in \eqref{440}, the definition of $\beta^{(s)}$ in \eqref{ac09}, and recalling that $l^{(\eps)}=k^{(V(0;\eps))}$, it follows that for any $\eps'\in[0,\iy)$, $\beta_{\eps'}=\beta^{(V(0;\eps'))}$. From the continuity of $\eps\mapsto V(\cdot;\eps)$, see Theorem \ref{thm_45}, and the second inequality in \eqref{435}, we get the first inequality in \eqref{460}.

The second inequality in \eqref{460} is trivial and follows since $\beta_{\eps+\delta}\le\hat\beta_{\eps+\delta}$. We now turn to proving the last inequality in \eqref{460}. 
Fix $\eps\in[0,\iy)$. From Theorem \ref{thm_45} and Lemma \ref{lem_43} we have
\begin{align}\notag
	\lim_{\delta\to 0}V(\cdot;\eps+\delta)=V(\cdot;\eps)\quad\text{and}\quad \lim_{\delta\to 0}l^{(\eps+\delta)}(\cdot)=l^{(\eps)}(\cdot),
\end{align}
uniformly on $[0,b]$. 
Together with \eqref{467} applied to $\eps'=\eps+\delta$, one has $V(x;\eps+\delta)=l^{(\eps+\delta)}(x)$, $x\in[0,\hat\beta_{\eps+\delta})$. Taking $\delta\to 0$, we get that $V(x;\eps)=l^{(\eps)}(x)$ for $x\in[0,\al_\eps)$, where $\al_\eps:=\limsup_{\delta\to 0}\hat\beta_{\eps+\delta}$. Since both functions are continuous, the equality holds for $\al_\eps$ as well. By the definition of $\hat\beta_\eps$, see \eqref{441}, we get that \eqref{460} holds.



\hfill$\Box$

\section*{Acknowledgments}
The author is greatful to an anonymous associate editor and to two
anonymous referees for their careful reading of the paper and for their suggestions that improved the presentation of the paper.

\footnotesize
\bibliographystyle{abbrv} 
\bibliography{refs} 

\end{document}